\documentclass[onefignum,onetabnum]{siamart251216}

\usepackage[utf8]{inputenc}

\usepackage{amssymb,amsfonts}
\usepackage{mathtools}
\usepackage{empheq}
\usepackage{extarrows}
\usepackage{dsfont}
\usepackage{stackengine}
\usepackage{hyperref}

\usepackage{graphicx}
\usepackage{float}
\usepackage{caption}
\usepackage{subcaption}
\usepackage{booktabs}
\usepackage{array}
\usepackage{multirow}
\usepackage{colortbl}      
\usepackage{placeins}

\usepackage{enumitem}
\usepackage{verbatim}
\usepackage{multicol}

\usepackage{algpseudocode}

\usepackage{tikz}
\usetikzlibrary{shapes,arrows}

\AtBeginDocument{%
  \setlength{\abovedisplayskip}{8pt plus 2pt minus 1pt}%
  \setlength{\belowdisplayskip}{8pt plus 2pt minus 1pt}%
  \setlength{\abovedisplayshortskip}{4pt plus 1pt minus 1pt}%
  \setlength{\belowdisplayshortskip}{4pt plus 1pt minus 1pt}%
}
\newsiamremark{remark}{Remark}
\crefname{remark}{Remark}{Remarks}

\newcommand{\be}{\begin{enumerate}}
\newcommand{\Ee}{\end{enumerate}}
\newcommand{\bc}{\begin{center}}
\newcommand{\ec}{\end{center}}

\newcommand{\bx}{\mathbf{x}}
\newcommand{\by}{\mathbf{y}}
\newcommand{\bu}{\mathbf{u}}
\newcommand{\bff}{\mathbf{f}}
\newcommand{\bg}{\mathbf{g}}

\newcommand{\R}{\mathbb{R}}
\newcommand{\cJ}{\mathcal{J}}

\headers{Hessian-augmented supervised learning for HJB PDEs}{G\'omez-Aedo, Azmi, Huang, Kalise, and Kunisch}

\title{Hessian-augmented Supervised Learning for Hamilton-Jacobi-Bellman PDEs}

\author{%
  Mat\'ias G\'omez-Aedo\thanks{Department of Mathematics, Imperial College
    London, South Kensington Campus, SW7 2AZ London, UK
    (\email{m.gomez-aedo22@imperial.ac.uk},
     \email{yuyang.huang21@imperial.ac.uk},
     \email{dkaliseb@imperial.ac.uk}).}
  \and Behzad Azmi\thanks{Department of Mathematics and Statistics,
    University of Konstanz, Germany
    (\email{behzad.azmi@uni-konstanz.de}).}
  \and Yuyang Huang\footnotemark[1]
  \and Dante Kalise\footnotemark[1]
  \and Karl Kunisch\thanks{RICAM and Institute of Mathematics and
    Scientific Computing, University of Graz, Austria
    (\email{karl.kunisch@uni-graz.at}).}}

\begin{document}
\maketitle
\begin{abstract}
A data-driven method is developed for approximating value functions in deterministic optimal control problems with nonlinear control-affine dynamics. The Pontryagin Maximum Principle optimality system is solved from multiple initial conditions to generate training data consisting of values, gradients, and Hessians of the value function, where Hessian information is obtained from a matrix Riccati equation along optimal trajectories. These quantities augment a weighted least-squares regression over sparse polynomial bases on hyperbolic cross index sets, with gradients and Hessians contributing additional linear equations per sample and substantially reducing sample complexity compared to value-only regression. Feedback laws are recovered analytically from the learned value function. In high dimensions, a partial Hessian strategy controls the cost of data generation. The approach is validated on problems of increasing state dimension, where second-order data augmentation is shown to improve approximation accuracy and closed-loop performance, with up to an order-of-magnitude reduction in the number of training samples required relative to lower-order methods.
\end{abstract}



\section{Introduction}\label{sec:int}
The synthesis of optimal feedback control laws is a fundamental challenge in control engineering and applied mathematics, arising across domains including aerospace guidance, robotic motion planning, mathematical finance, and the stabilisation of PDE-governed systems. The standard approach begins by casting a dynamic optimisation problem
\begin{align}\label{eq:optcostinf_intro}
	\underset{\mathbf{u}(\cdot)\in L^2([0,\infty);\mathbb{R}^m)}{\min}\; 
	\cJ(\mathbf{u};\mathbf{x}) := \int^{\infty}_{0} 
	\ell(\mathbf{y}(s))+\beta\|\mathbf{u}(s)\|_2^2\,ds, \qquad \beta>0,
\end{align}
where $\mathbf{u}(s) \in \mathbb{R}^m$ is a time-dependent control signal and the state $\mathbf{y}(s)\in\mathbb{R}^n$ satisfies the control-affine nonlinear dynamics
\begin{align}\label{eq:state_intro}
	\frac{d}{ds}\mathbf{y}(s) = \mathbf{f}(\mathbf{y}(s)) 
	+ \mathbf{g}(\mathbf{y}(s))\mathbf{u}(s),
\end{align}
with initial condition $\by(0)=\bx$. The control objective is encoded in the state penalty $\ell(\mathbf{y})$, and the dimensions $(n,m)$ can be large: a quadrotor model involves up to $n=12$ states, agent-based systems (swarm robotics, financial portfolios) reach hundreds, and fluid flow control via Navier--Stokes constitutes an infinite-dimensional state space. Our concern is with real-time, robust feedback laws of the form $\bu=\bu(\by(s))$, which depend on the current state rather than on time and initial condition alone. By the dynamic programming principle (DPP), the optimal feedback map is given by
\begin{align}\label{eq:optf_intro}
	\bu(\bx) = -\frac{1}{2\beta}\bg(\bx)^{\top}\nabla V(\bx),
\end{align}
where $V: \R^n\to\R$  is the \emph{value function} (optimal cost-to-go), defined as 
\begin{align}\label{eq:value_intro}
	V(\bx) := \underset{\mathbf{u}(\cdot)}{\inf}\;\cJ(\mathbf{u}(\cdot);\mathbf{x}),
\end{align}
which solves the Hamilton--Jacobi--Bellman (HJB) equation
\begin{align}\label{eq:thjb_intro}
	\underset{\bu\in \R^m}{\min}\bigl\{ 
	(\bff(\bx)+\bg(\bx)\bu)^\top \nabla V(\bx) 
	+ \ell(\bx)+\beta\|\bu\|_2^2 \bigr\} = 0.
\end{align}

The feedback law \eqref{eq:optf_intro} is globally optimal and compensates disturbances along trajectories without recomputing the control action. Its chief drawback is computational: the HJB equation is a fully nonlinear, non-variational PDE posed over the $n$-dimensional state space, and classical grid-based discretisations scale exponentially in $n$, limiting their scope to problems with $n \lesssim 5$.

Circumventing the curse of dimensionality has motivated a broad range of non-standard numerical strategies. Causality-free characteristics methods exploit the connection between optimal trajectories and HJB characteristics to recover feedback laws without state-space discretisation \cite{KW1,KW2,Yegorov,SUB}, while Lax--Hopf formulas combined with convex optimisation enable pointwise value function evaluation in high dimensions \cite{DO}. Sparse polynomial approximations from Galerkin and policy-iteration schemes have been applied to HJB-based stabilization of semilinear parabolic systems \cite{KK2018}. Al'brekht's power-series method, reformulated via Kronecker-product tensor structure, yields scalable polynomial approximations of the HJB equation for polynomial control-affine systems, reaching state dimensions exceeding one thousand \cite{CK1,CK2}. Low-rank tensor decompositions exploit separable structure in the value function \cite{DKK,DKS}, and separable neural-network representations have been recently pursued along similar lines \cite{SSKG}. Neural networks trained by residual minimisation over monotone finite-difference discretisations \cite{ETM} and direct training over the feedback map \cite{KW21} offer complementary perspectives. Deep learning approaches, including physics-informed and Deep Galerkin methods \cite{Sirignano}, backward SDE formulations \cite{HJE}, and symplectic networks applied to the PMP system \cite{MZDK}, have extended this to dimensions exceeding one hundred.

This paper adopts a \emph{derivative-augmented} supervised learning framework in which a polynomial model $\tilde V(\mathbf{x}) = \sum_{k=1}^q \theta_k \Phi_{i_k}(\mathbf{x})$ is fitted to a labelled dataset carrying values, gradients, and Hessians of $V$,
\begin{equation}\label{dataset_intro}
\mathcal{D} = \bigl\{(\mathbf{x}^j,\, V(\mathbf{x}^j),\,\nabla V(\mathbf{x}^j),\,
\nabla^2 V(\mathbf{x}^j))\bigr\}_{j=1}^N,
\end{equation}
by minimising
\begin{multline}\label{eq:intro_loss}
\min_{\theta \in \mathbb{R}^q}\;\frac{1}{N}\sum_{j=1}^N\Bigl[
 \bigl|\tilde V(\mathbf{x}^j) - V(\mathbf{x}^j)\bigr|^2
+ \gamma_1 \sum_{m=1}^{n}\bigl|\partial_{x_m}\tilde V(\mathbf{x}^j) 
  - \partial_{x_m}V(\mathbf{x}^j)\bigr|^2 \\
+ \gamma_2 \sum_{1 \le m \le \ell \le n}
  \bigl|\partial^2_{x_m x_\ell}\tilde V(\mathbf{x}^j) 
  - \partial^2_{x_m x_\ell}V(\mathbf{x}^j)\bigr|^2\Bigr],
\end{multline}
where $\gamma_1,\gamma_2\ge 0$ weight gradient and Hessian contributions, and $\{\Phi_{i_k}\}$ is a polynomial basis over a finite multi-index set. For control-affine dynamics with quadratic control cost, the feedback law is recovered analytically from \eqref{eq:optf_intro}.

The dataset is generated by solving the first-order necessary optimality conditions of the Pontryagin Maximum Principle (PMP) \cite{Pontryagin} from a collection of initial conditions. Each PMP solve is a two-point boundary value problem (TPBVP) that can be treated independently of any spatial grid, avoiding the exponential scaling of state-space discretisation. When the PMP conditions are also sufficient, the value function can be evaluated by following optimal trajectories; the adjoint variable then coincides with $\nabla V$, and $\nabla^2 V$ is recoverable from a matrix  Riccati equation integrated alongside the same ODE system.

Derivative-augmented approximation dates back to Hermite interpolation; its extension to high-dimensional polynomial approximation from gradient measurements was studied in \cite{AS}, and its use in neural network training is known as Sobolev training \cite{COJSP}. In the HJB setting, gradient-augmented regression was introduced in \cite{AKK,NGK} and subsequently extended to Hermite kernel surrogates \cite{EH}, gradient-enhanced tensor train cross approximation \cite{DKS}, RKHS approximations driven by HJB verification conditions \cite{EAH}, and supervised learning from State-dependent Riccati Equation solves \cite{ABK}. Hessian augmentation has appeared in the context of the Mortensen observer \cite{BKS}, where second-order data improves value function approximation. The present work systematically develops this direction for optimal control.

\begin{figure}[t]
\centering
\begin{minipage}[t]{0.32\textwidth}%
    \includegraphics[width=\linewidth]{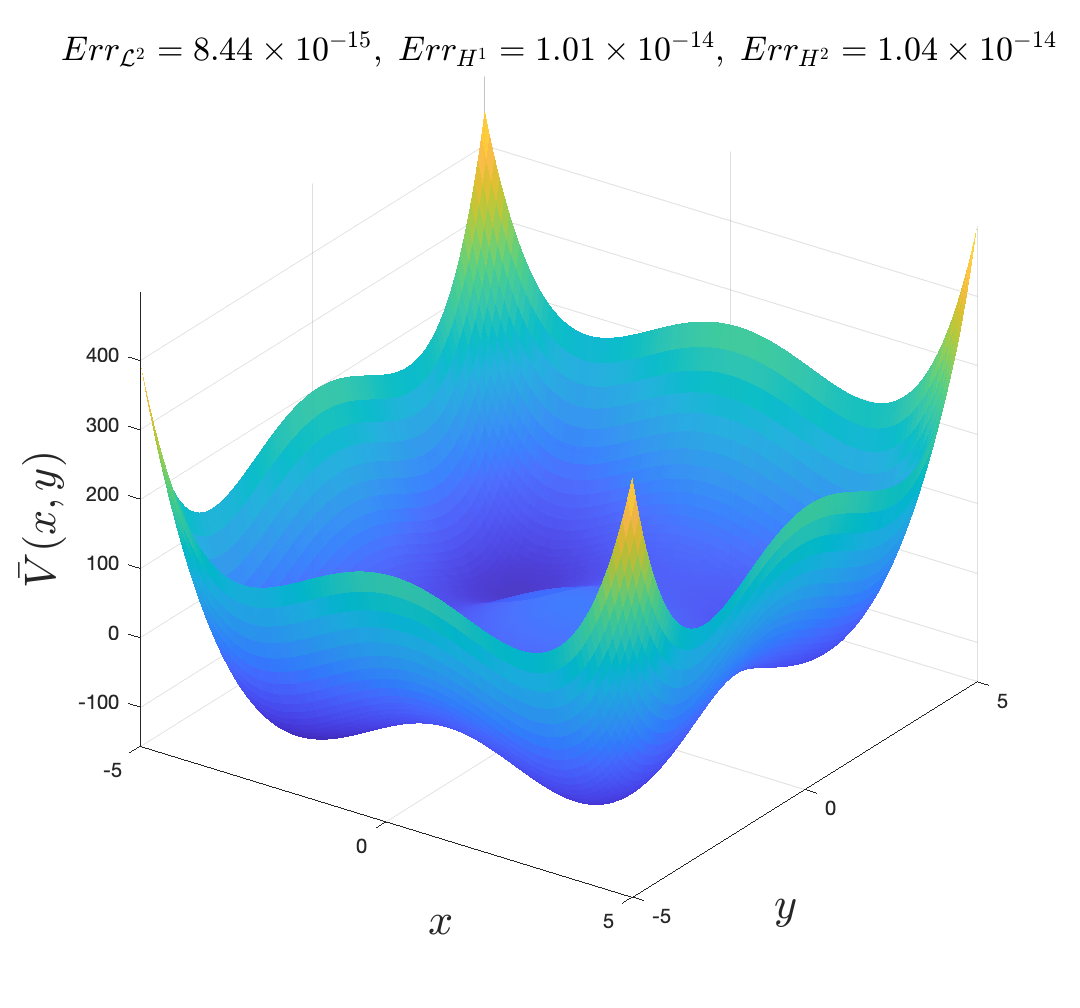}
    \subcaption*{(i) $N = 21$, zeroth-order}
\end{minipage}\hfill
\begin{minipage}[t]{0.32\textwidth}%
    \includegraphics[width=\linewidth]{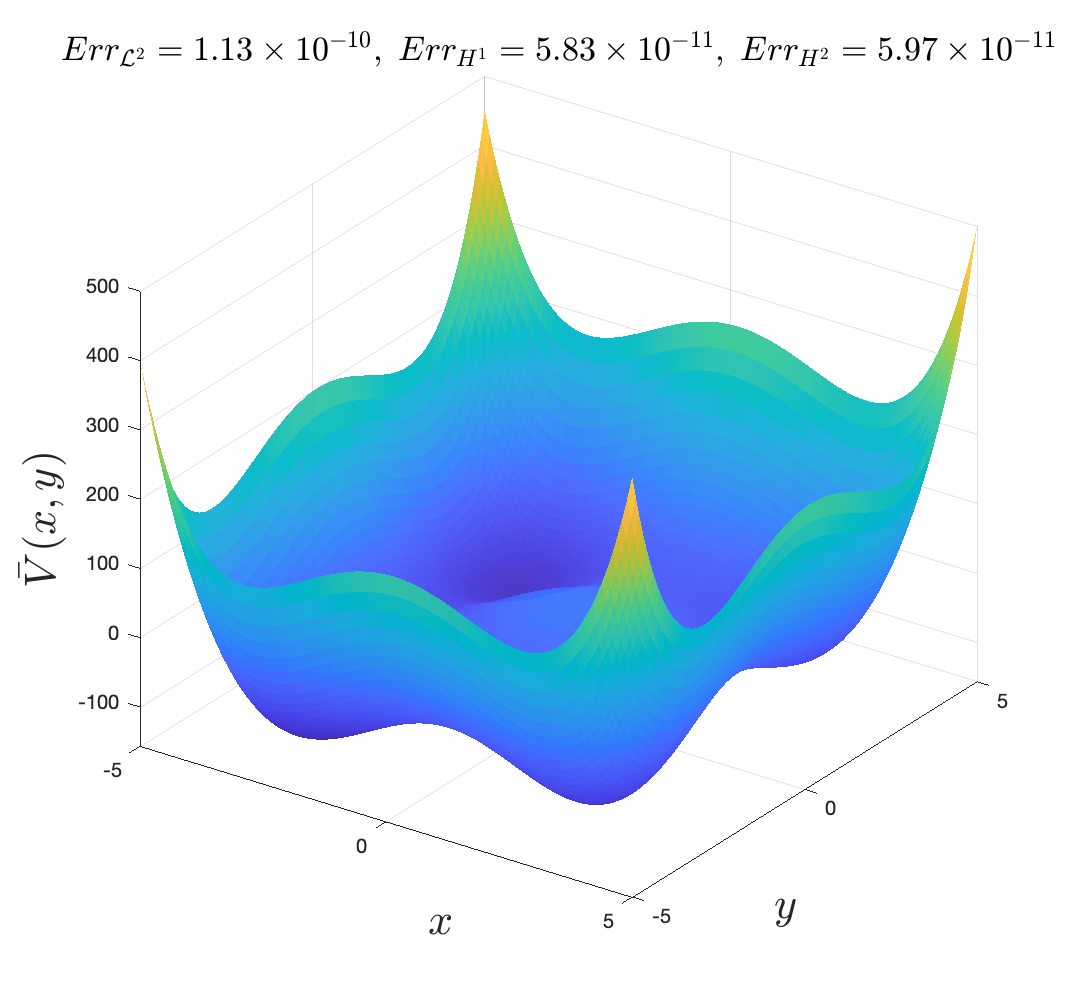}
    \subcaption*{(ii) $N = 7$, first-order}
\end{minipage}\hfill
\begin{minipage}[t]{0.32\textwidth}%
    \includegraphics[width=\linewidth]{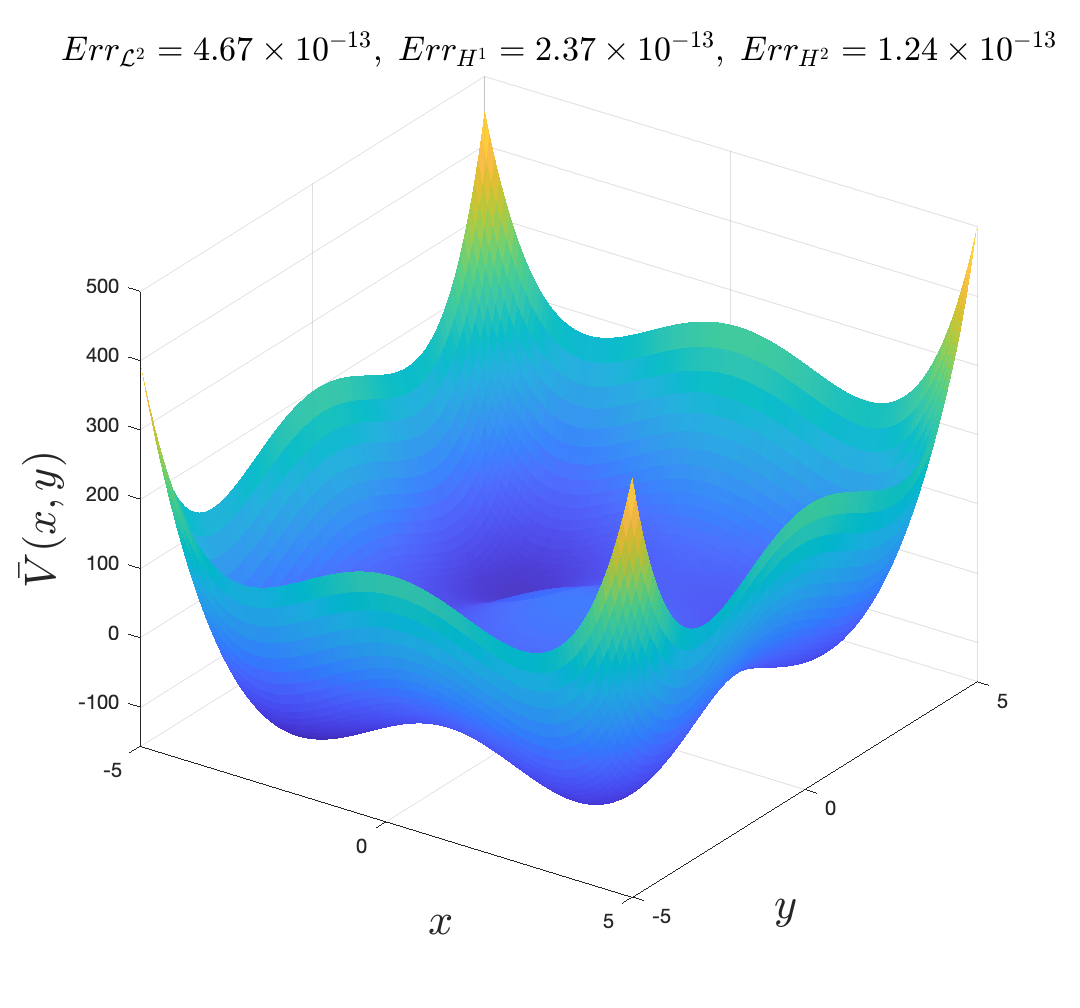}
    \subcaption*{(iii) $N = 4$, second-order}
\end{minipage}
\caption{Approximation of the Styblinski--Tang function 
    $f_{ST}(x,y) = \tfrac{1}{2}(x^4-16x^2+5x+y^4-16y^2+5y)$ using zeroth-, 
    first-, and second-order augmentation, at the minimum training set size 
    achieving a prescribed accuracy on a validation set.}
\label{fig:ST}
\end{figure}

The main motivation for incorporating Hessian data is reducing sample complexity. Gradient data contributes $n$ equations per sample and Hessian data a further $n(n+1)/2$, giving a total equation count of $M = N(1+n+n(n+1)/2)$ in \eqref{eq:intro_loss}, a multiplicative gain that can be decisive when each PMP data point requires solving a $2n$-dimensional TPBVP. Figure~\ref{fig:ST} illustrates this concretely for polynomial regression of a standard test function: the minimum sample count achieving a prescribed accuracy drops from $N=21$ 
(values only) to $N=7$ (gradients) to $N=4$ (Hessians).
The contributions of the paper are threefold. We extend the PMP-based data generation approach of \cite{AKK} to second-order information  via Riccati equations, and introduce a \emph{partial Hessian} strategy in which second-order data is included for a tunable fraction $\rho\in(0,1]$ of training  points. We study the interaction between second-order augmentation and hyperbolic cross polynomial bases enriched with low-degree total-degree terms. Finally, we conduct a systematic benchmark across problems of increasing state 
dimension, measuring accuracy in discrete Sobolev norms, HJB residual, and closed-loop performance.

\noindent The paper is organised as follows. Section~\ref{sec:ocp} formulates the optimal control problem and derives the PMP-based data generation pipeline. Section~\ref{sec:method} describes the regression framework, basis selection, partial Hessian formulation, and error metrics. Section~\ref{sec:numerics} reports numerical experiments on the Van der Pol oscillator, a satellite attitude control problem, and a PDE-constrained problem arising from the Allen--Cahn equation.

\FloatBarrier

\section{Optimal control, HJB equation, and PMP-based data generation}
\label{sec:ocp}
This section establishes the analytical framework that underpins our data generation pipeline. We reformulate the optimal control problem on a finite horizon and characterize the value function via the HJB equation. We then state the PMP and make precise the correspondence between the two approaches, which is the structural assumption for extracting value, gradient, and Hessian data from PMP solves.

\subsection{Finite-horizon formulation and the value function}
We consider the finite-horizon optimal control problem
\begin{align}\label{eq:optcost}
	\underset{\mathbf{u}(\cdot)\in \mathcal{U}}{\min}\;
	J(\mathbf{u};t,\mathbf{x}) := \int^{T}_{t} 
	\ell(\mathbf{y}(s))+\beta\|\mathbf{u}(s)\|^2\,ds,\quad t\in[0,T),\;\beta>0,
\end{align}
with $\mathcal{U}:=L^2([t,T];\mathbb{R}^m)$, subject to the control-affine dynamics \eqref{eq:state_intro} and initial condition $\by(t)=\bx$. We assume throughout that $\ell$, $\mathbf{f}$, and $\mathbf{g}$ are continuously differentiable. The value function $V:[0,T]\times\mathbb{R}^n\to\mathbb{R}$ is defined as the infimum of the cost over admissible controls,
\begin{equation}\label{eq:value}
	V(t,\mathbf{x}) := \underset{\mathbf{u}(\cdot)\in\mathcal{U}} {\inf}
	\bigl\{J(\mathbf{u};t,\mathbf{x}) \text{ subject to } \eqref{eq:state_intro}\;\text{and}\;\by(t)=\bx\bigr\},
\end{equation}
encodes the optimal cost-to-go from every state $\mathbf{x}$ at every time $t$. The DPP \cite[Chapter 1]{BCD} implies that $V$ is the unique viscosity solution of the HJB equation
\begin{equation}\label{HJBe}
	\begin{cases}
		\partial_t V(t,\mathbf{x}) 
		+ \nabla_{\bx} V(t,\mathbf{x})^{\top}\mathbf{f}(\mathbf{x})
		- \dfrac{1}{4\beta}\|\mathbf{g}(\mathbf{x})^{\top}\nabla_{\bx} V(t,\mathbf{x})\|_2^2
		+ \ell(\mathbf{x}) = 0,\\[6pt]
		V(T,\mathbf{x}) = 0,
	\end{cases}
\end{equation}
and the optimal feedback law is recovered from its spatial gradient as
\begin{equation}\label{eq:feedhjb}
	\mathbf{u}^*(t,\mathbf{x}) 
	= -\frac{1}{2\beta}\mathbf{g}(\mathbf{x})^{\top}\nabla_{\bx} V(t,\mathbf{x}).
\end{equation}
This feedback law is time-dependent, unlike the stationary policies targeted in Section~\ref{sec:int}. For $T$ sufficiently large, however, $V(0,\mathbf{x})$ accurately approximates the infinite-horizon value $V(\mathbf{x})$ in
\eqref{eq:value_intro} \cite{GR08}, so the map $\mathbf{x}\mapsto\mathbf{u}^*(0,\mathbf{x})$ provides the static feedback law of the form \eqref{eq:optf_intro} we seek. Direct numerical solution of \eqref{HJBe} to obtain $V(0,\cdot)$ is however intractable in high dimensions due to the exponential scaling of any state-space discretisation, which motivates the PMP-based data generation strategy developed next.

\subsection{Pontryagin's Maximum Principle and PMP-based data generation}
Our strategy is to generate labelled data for $V(0,\cdot)$ by solving the first-order optimality conditions of \eqref{eq:optcost} via PMP, and to fit a global approximation $\tilde{V}\approx V(0,\cdot)$ from which the feedback law is recovered using \eqref{eq:feedhjb}. For the problem \eqref{eq:optcost}, define the PMP Hamiltonian
\begin{equation}\label{eq:Hpmp}
    H(\mathbf{y},\mathbf{p},\mathbf{u})
    := \ell(\mathbf{y}) + \beta\|\mathbf{u}\|_2^2
    + \mathbf{p}^{\top}\bigl(\mathbf{f}(\mathbf{y})+\mathbf{g}(\mathbf{y})\mathbf{u}\bigr).
\end{equation}
The PMP characterises optimal triples $(\mathbf{y}^*(s),\mathbf{p}^*(s),\mathbf{u}^*(s))$ via the canonical system
\begin{equation*}
    \frac{d}{ds}\mathbf{y}^*(s) = \nabla_{\mathbf{p}} H(\mathbf{y}^*(s),\mathbf{p}^*(s),\mathbf{u}^*(s)), \qquad -\frac{d}{ds}\mathbf{p}^*(s) = \nabla_{\mathbf{y}} H(\mathbf{y}^*(s),\mathbf{p}^*(s),\mathbf{u}^*(s)),
\end{equation*}
together with the pointwise optimality condition $\mathbf{u}^*(s) = \underset{\bu}{\operatorname{argmin}}  \,H(\mathbf{y}^*(s),\mathbf{p}^*(s),\mathbf{u})$. Since the control enters the dynamics affinely and the control cost is quadratic, this yields explicitly
\begin{equation}\label{eq:optcon}
    \mathbf{u}^*(s) = -\frac{1}{2\beta}\mathbf{g}(\mathbf{y}^*(s))^{\top}\mathbf{p}^*(s),
\end{equation}
and the reduced Hamiltonian
\begin{equation}\label{eq:Hreduced}
    \mathcal{H}(\mathbf{y},\mathbf{p})
    := \min_{\mathbf{u}\in\mathbb{R}^m} H(\mathbf{y},\mathbf{p},\mathbf{u})
    = \ell(\mathbf{y})
    + \mathbf{p}^{\top}\mathbf{f}(\mathbf{y})
    - \frac{1}{4\beta}\|\mathbf{g}(\mathbf{y})^{\top}\mathbf{p}\|_2^2.
\end{equation}
Substituting \eqref{eq:optcon}-\eqref{eq:Hreduced} into the canonical equations gives the $2n$-dimensional TPBVP
\begin{align}\label{eq:tpbvp}
    \left\{\begin{aligned}
        \frac{d}{ds}\mathbf{y}^*(s) &= \nabla_{\mathbf{p}}\mathcal{H}(\mathbf{y}^*(s),\mathbf{p}^*(s)),
        & \mathbf{y}^*(t) &= \mathbf{x},\\[4pt]
        -\frac{d}{ds}\mathbf{p}^*(s) &= \nabla_{\mathbf{y}}\mathcal{H}(\mathbf{y}^*(s),\mathbf{p}^*(s)),
        & \mathbf{p}^*(T) &= 0.
    \end{aligned}\right.
\end{align}
The dataset \eqref{dataset_intro} is constructed by sampling initial conditions $(t, \mathbf{x})$ independently and solving \eqref{eq:tpbvp}, without enforcing the causal structure of the HJB equation.

\noindent The PMP system \eqref{eq:tpbvp} and the HJB equation \eqref{HJBe} are linked by the method of characteristics \cite{SUB}: the optimal state-adjoint system coincides with the characteristic equations of the HJB, and solving it from each initial condition $\mathbf{x}$ propagates information about $V$ along the corresponding optimal trajectory. This connection is fundamental to all three levels of data extraction formalised below.

\begin{proposition}[Zeroth- and first-order relations]\label{prop:pmp_dpp}
Suppose $V$ is a classical solution of \eqref{HJBe} of class $\mathcal{C}^1$,
and that the PMP system \eqref{eq:tpbvp} admits a unique solution
$(\mathbf{y}^*,\mathbf{p}^*,\mathbf{u}^*)$ from initial condition
$\mathbf{y}^*(t)=\mathbf{x}$. Then, for all $s\in[t,T]$:
\begin{enumerate}
    \item[\textnormal{(i)}]  $V(s,\mathbf{y}^*(s)) = J(\mathbf{u}^*;s,\mathbf{y}^*(s))$, where $\mathbf{u}^*$ is restricted to $[s,T]$,
    \item[\textnormal{(ii)}] $\nabla_{\mathbf{x}} V(s,\mathbf{y}^*(s))=\mathbf{p}^*(s)$.
\end{enumerate}
\end{proposition}
While Proposition~\ref{prop:pmp_dpp} extracts values and gradients of $V$ from the PMP solution, second-order information is accessible through the same trajectory via a matrix Riccati equation for $\nabla^2_{\mathbf{x}}V(s,\mathbf{y}^*(s))$.
\begin{theorem}[Second-order relations \cite{CV,CV1983,CFS}]\label{thm:riccati}
Assume that $V(s,\cdot)$ is twice Fr\'echet differentiable at $\mathbf{y}^*(s)$
for all $s\in[t,T]$, and that the reduced Hamiltonian $\mathcal{H}$ defined
in \eqref{eq:Hreduced} is sufficiently smooth in a neighbourhood of
$\{(\mathbf{y}^*(s),\mathbf{p}^*(s)):s\in[t,T]\}$.
Then $R(s):=\nabla^2_{\mathbf{x}}V(s,\mathbf{y}^*(s))$ satisfies the
matrix Riccati equation
\begin{equation}\label{eq:riccati}
    \begin{cases}
    \dot{R}(s) + \mathcal{H}_{\mathbf{py}}[s]\,R(s) + R(s)\,\mathcal{H}_{\mathbf{yp}}[s]
    + R(s)\,\mathcal{H}_{\mathbf{pp}}[s]\,R(s) + \mathcal{H}_{\mathbf{yy}}[s] = 0,\\
    R(T) = 0,
    \end{cases}
\end{equation}
where $\mathcal{H}_{\mathbf{py}}[s] :=
\nabla^2_{\mathbf{py}}\mathcal{H}(\mathbf{y}^*(s),\mathbf{p}^*(s))$,
and similarly for $\mathcal{H}_{\mathbf{yp}}[s]$, $\mathcal{H}_{\mathbf{pp}}[s]$,
$\mathcal{H}_{\mathbf{yy}}[s]$. 
\end{theorem}
Identity (i) is the standard verification theorem; (ii) is the costate identification, rigorously established via the method of characteristics; the second-order Riccati identity \eqref{eq:riccati} was first derived by Clarke and Vinter \cite{CV,CV1983} and further developed by Cannarsa, Frankowska, and collaborators \cite{CF,CFS,CF92,CF96}. These references cover the problem in greater generality, addressing the non-smooth case via Fréchet subdifferentials. Our formulation suffices for the smooth control problems considered in Section~\ref{sec:numerics}.

\begin{remark}[Conjugate points]\label{rem:conjugate}
	The Riccati equation \eqref{eq:riccati} develops a finite-time singularity at conjugate points \cite{CF96}. Geometrically, conjugate points mark the merging of distinct optimal trajectories and the loss of differentiability of $V$. In such cases, the Hessian data and the corresponding samples should be excluded from the training dataset. For the control problems considered in Section~\ref{sec:numerics}, the absence of conjugate points on the sampled trajectories is verified numerically by monitoring the condition number of $R(t)$ throughout integration.
\end{remark}

\FloatBarrier

\section{Second-order, data-driven approximation of the value function}
\label{sec:method}
This section describes the computational framework for approximating the value function from the augmented dataset \eqref{dataset_intro} and recovering the associated feedback law. We detail the polynomial approximation architecture, the derivative-augmented regression problem and its partial Hessian variant, and the error metrics used for assessment. For notational brevity, we fix 
$t = 0$ throughout and write $V(\mathbf{x}) := V(0,\mathbf{x})$, 
$\nabla V(\mathbf{x}) := \nabla_{\mathbf{x}} V(0,\mathbf{x})$, 
and $\nabla^2 V(\mathbf{x}) := \nabla^2_{\mathbf{x}} V(0,\mathbf{x})$.

\subsection{Data generation}
\label{sec:augmented_data}
A single solve of the PMP system \eqref{eq:tpbvp} from initial condition $\mathbf{x}^j$, combined with the backward integration of the Riccati equation \eqref{eq:riccati} along the resulting optimal trajectory, yields the triple $(V(\mathbf{x}^j), \nabla V(\mathbf{x}^j), \nabla^2 V(\mathbf{x}^j))$ as detailed in Proposition~\ref{prop:pmp_dpp} and Theorem~\ref{thm:riccati}. In practice, the TPBVP \eqref{eq:tpbvp} is solved with a dedicated boundary value problem solver (e.g.\ \texttt{bvp5c} in \textsc{Matlab}), and the Riccati equation \eqref{eq:riccati} is integrated backward with a standard ODE solver (e.g.\ \texttt{ode45} in \textsc{Matlab}) over the same time grid. Hessian samples are discarded whenever the condition number of $R(0)$ indicates proximity to a conjugate point. Repeating for $N$ initial conditions $\{\mathbf{x}^j\}_{j=1}^N$ sampled from a domain $\Omega \subset \mathbb{R}^n$ assembles the augmented dataset
\begin{equation}\label{dataset2}
    \mathcal{D} = \bigl\{
    (\mathbf{x}^j,\, V(\mathbf{x}^j),\, \nabla V(\mathbf{x}^j),\, 
    \nabla^2 V(\mathbf{x}^j))
    \bigr\}_{j=1}^N.
\end{equation}
For each test, the initial conditions $\{\mathbf{x}^j\}_{j=1}^N$ are sampled from an $n$-dimensional hyperrectangle $\Omega \subset \mathbb{R}^n$ using a scrambled Sobol sequence, distributing the points quasi-uniformly over $\Omega$. The regression itself imposes no requirement on the sampling scheme.

\subsection{Polynomial models for data-driven approximation}
\label{sec:poly_models}

We approximate the value function $\tilde{V}:\mathbb{R}^n\to\mathbb{R}$ on a bounded domain $\Omega\subset\mathbb{R}^n$ by a linear combination of multivariate polynomial basis functions. For $\Omega=(-1,1)^n$, let $\{\phi_i\}_{i=0}^\infty$ be a one-dimensional orthonormal Legendre basis\footnote{The choice of Legendre polynomials throughout this paper is
made for simplicity and is not essential: any polynomial basis, orthogonal
or monomial, can be used in the regression framework. Orthonormal bases are
preferable for the conditioning of the design matrix.} of $L^2(-1,1)$. The tensor-product basis $\{\Phi_\mathbf{i}\}_{\mathbf{i}\in\mathbb{N}_0^n}$ of $L^2(\Omega)$ is
\begin{equation}\label{eq:polbasis}
    \Phi_\mathbf{i}(\mathbf{x}) := \prod_{j=1}^n \phi_{i_j}(x_j), \qquad 
    \mathbf{i}=(i_1,\ldots,i_n)\in\mathbb{N}_0^n.
\end{equation}
Given a finite multi-index set $\mathfrak{I}\subset\mathbb{N}_0^n$ with 
$|\mathfrak{I}|=:q<\infty$, the approximation takes the form
\begin{equation}\label{eq:approxValue}
    \tilde{V}(\mathbf{x}) = \sum_{\mathbf{i}\in\mathfrak{I}} 
    \theta_\mathbf{i}\,\Phi_\mathbf{i}(\mathbf{x}), \qquad 
    \theta_{\mathfrak{I}} := \{\theta_\mathbf{i}\}_{\mathbf{i}\in\mathfrak{I}} 
    \in \mathbb{R}^q,
\end{equation}
from which the feedback law is recovered analytically via \eqref{eq:feedhjb}. The central design choice is the index set $\mathfrak{I}$: it must remain tractable as $n$ grows while capturing the relevant structure of $V$. We consider two constructions.
\noindent\textit{Total degree.} For a given degree $s$, the multi-indices are those whose components sum to at most $s$:
\begin{equation}\label{eq:totaldegree}
    \mathfrak{I}_{\mathrm{TD}}(s) = \{\mathbf{i}\in\mathbb{N}_0^n:
    \|\mathbf{i}\|_1\leq s\}, \qquad 
    |\mathfrak{I}_{\mathrm{TD}}(s)| = \sum_{j=0}^s\binom{n+j-1}{j}.
\end{equation}
\noindent\textit{Hyperbolic cross.} The indices satisfy a product constraint, yielding sub-exponential cardinality \cite{ABW_book}:
\begin{equation}\label{eq:hypcross}
    \mathfrak{I}_{\mathrm{HC}}(s) = \Bigl\{\mathbf{i}\in\mathbb{N}_0^n:
    \prod_{j=1}^n(i_j+1)\leq s+1\Bigr\}, \qquad 
    |\mathfrak{I}_{\mathrm{HC}}(s)| \leq 
    \min\bigl\{2s^3 4^n,\, e^2 s^{2+\log_2 n}\bigr\}.
\end{equation}
The hyperbolic cross is nested: increasing $s$ adds multi-indices without removing existing ones. We adopt it as our primary index set.
In practice, we enrich the hyperbolic cross with the total degree set at a low degree $s_{\mathrm{td}}$:
\begin{equation}\label{eq:enriched_index}
    \mathfrak{I} = \mathfrak{I}_{\mathrm{HC}}(s)\cup
    \mathfrak{I}_{\mathrm{TD}}(s_{\mathrm{td}}).
\end{equation}
The hyperbolic cross at moderate $s$ favours multi-indices concentrated in few dimensions but may miss low-order cross-terms involving several coordinates simultaneously. Such terms are essential when the value function contains genuine multivariate coupling of low individual degree. For instance, the monomial $x_1 x_2 x_3$ is present in the total-degree basis with $s_{\mathrm{td}}=3$, but is excluded from the hyperbolic cross at any $s<7$. Enrichment at $s_{\mathrm{td}}=5$--$7$ adds a polynomially-scaling number of terms via \eqref{eq:totaldegree} without compromising the dimensional scaling of \eqref{eq:hypcross}. The interaction between this enrichment and second-order augmentation is examined in Section~\ref{sec:numerics}.

\begin{remark}[Domain rescaling]\label{rem:rescaling}
When $\Omega=\prod_{j=1}^n[-b_j,b_j]$, each coordinate is rescaled a $\tilde{x}_j = x_j/a_j$. Setting $a_j = b_j$ maps $\Omega$ exactly onto
$(-1,1)^n$, but places sample points near $\pm 1$ where high-degree Legendre
polynomials attain their largest values, potentially degrading the conditioning
of the regression matrix. Since sample points are drawn from $\Omega$, setting
$a_j = b_j(1+\mu)$ for a small margin $\mu>0$ maps all samples into
$[-1/(1+\mu),\,1/(1+\mu)]^n$, keeping them bounded away from $\pm 1$.
\end{remark}

\subsection{Regression with derivative-augmented data}
\label{sec:regression}
Given the augmented dataset $\mathcal{D}$ from \eqref{dataset2}, we construct 
a weighted least-squares system incorporating zeroth-, first-, and second-order 
information. The coefficient vector $\theta_{\mathfrak{I}}\in\mathbb{R}^q$ in 
\eqref{eq:approxValue} is determined by minimising a quadratic objective that 
penalises misfit in function values, gradients, and Hessians simultaneously.

For each sample point $\mathbf{x}^j$, we introduce the abbreviations
\[
V^j := V(\mathbf{x}^j), 
\qquad
\partial_{x_m} V^j := \partial_{x_m} V(\mathbf{x}^j),
\qquad
\partial^2_{x_m x_\ell} V^j := \partial^2_{x_m x_\ell} V(\mathbf{x}^j).
\]
Since the Hessian is symmetric, we retain only the distinct second-order entries corresponding to $1 \le m \le \ell \le n$. We fix an enumeration $\mathfrak{I}=\{\mathbf{i}_1,\ldots,\mathbf{i}_q\}$ and identify  $\theta_{\mathfrak{I}}$ with the vector $\theta=(\theta_{\mathbf{i}_1},\ldots \theta_{\mathbf{i}_q})^\top\in\mathbb{R}^q$. For the function values, we define 
$\mathbf{V}_0 \in \mathbb{R}^N$ and $\mathbf{A}_0 \in \mathbb{R}^{N \times q}$ by
\begin{equation}\label{eq:block0}
\mathbf{V}_0 := \tfrac{1}{\sqrt{N}}\bigl(V^j\bigr)_{j=1}^N,
\qquad
\mathbf{A}_0 := \tfrac{1}{\sqrt{N}}\bigl(\Phi_{\mathbf{i}_k}(\mathbf{x}^j)
\bigr)_{j,k=1}^{N,q}.
\end{equation}
For each coordinate direction $m=1,\dots,n$, the first-order blocks 
$\mathbf{V}_m \in \mathbb{R}^N$ and $\mathbf{A}_m \in \mathbb{R}^{N \times q}$ are
\begin{equation}\label{eq:block1}
\mathbf{V}_m := \tfrac{1}{\sqrt{N}}\bigl(\partial_{x_m}V^j
\bigr)_{j=1}^N,
\qquad
\mathbf{A}_m := \tfrac{1}{\sqrt{N}}\bigl(\partial_{x_m}\Phi_{\mathbf{i}_k}(\mathbf{x}^j)
\bigr)_{j,k=1}^{N,q}.
\end{equation}
Finally, for each pair $1\le m\le\ell\le n$, the second-order blocks 
$\mathbf{V}_{m,\ell} \in \mathbb{R}^N$ and 
$\mathbf{A}_{m,\ell} \in \mathbb{R}^{N \times q}$ are
\begin{equation}\label{eq:block2}
\mathbf{V}_{m,\ell} := \tfrac{1}{\sqrt{N}}\bigl(\partial^2_{x_mx_\ell}V^j
\bigr)_{j=1}^N,
\qquad
\mathbf{A}_{m,\ell} := \tfrac{1}{\sqrt{N}}\bigl(
\partial^2_{x_mx_\ell}\Phi_{\mathbf{i}_k}(\mathbf{x}^j)
\bigr)_{j,k=1}^{N,q}.
\end{equation}
The $1/\sqrt{N}$ normalisation ensures the objective scales as an empirical 
mean, independently of $N$. Stacking these blocks with weights 
$\gamma_1,\gamma_2\ge 0$,
\begin{equation}\label{eq:weighted_system}
\bar{\mathbf{A}}_{\gamma_1,\gamma_2}
:=
\left(\begin{array}{@{}r@{}l@{}}
& \mathbf{A}_0 \\
\sqrt{\gamma_1}\, & \mathbf{A}_1 \\
& \vdots \\
\sqrt{\gamma_1}\, & \mathbf{A}_n \\
\sqrt{\gamma_2}\, & \mathbf{A}_{1,1} \\
\sqrt{\gamma_2}\, & \mathbf{A}_{1,2} \\
& \vdots \\
\sqrt{\gamma_2}\, & \mathbf{A}_{n,n}
\end{array}\right),
\qquad
\bar{\mathbf{V}}_{\gamma_1,\gamma_2}
:=
\left(\begin{array}{@{}r@{}l@{}}
& \mathbf{V}_0 \\
\sqrt{\gamma_1}\, & \mathbf{V}_1 \\
& \vdots \\
\sqrt{\gamma_1}\, & \mathbf{V}_n \\
\sqrt{\gamma_2}\, & \mathbf{V}_{1,1} \\
\sqrt{\gamma_2}\, & \mathbf{V}_{1,2} \\
& \vdots \\
\sqrt{\gamma_2}\, & \mathbf{V}_{n,n}
\end{array}\right),
\end{equation}
where pairs $(m,\ell)$ follow the ordering of \eqref{eq:block2}, the 
derivative-augmented regression problem reads
\begin{equation}\label{eq:weighted_leastsquares_compact}
\min_{\theta\in\mathbb{R}^q}
\bigl\|\bar{\mathbf{A}}_{\gamma_1,\gamma_2}\theta
-\bar{\mathbf{V}}_{\gamma_1,\gamma_2}\bigr\|_2^2,
\end{equation}
with design matrix $\bar{\mathbf{A}}_{\gamma_1,\gamma_2}\in\mathbb{R}^{M\times q}$ 
and
\[
M := N\!\left(1+n+\frac{n(n+1)}{2}\right).
\]
Problem \eqref{eq:weighted_leastsquares_compact} is solved via QR factorisation. When the columns of 
$\bar{\mathbf{A}}_{\gamma_1,\gamma_2}$ are ill-conditioned, we consider the $\ell_2$-regularised formulation
\begin{equation}\label{eq:ridgeprob_reg}
\min_{\theta \in \mathbb{R}^q}
\bigl\|\bar{\mathbf{A}}_{\gamma_1,\gamma_2}\theta
-\bar{\mathbf{V}}_{\gamma_1,\gamma_2}\bigr\|_2^2
+\lambda\|\theta\|_2^2,
\qquad \lambda>0.
\end{equation}
Once $\theta^*$ is obtained from either \eqref{eq:weighted_leastsquares_compact} 
or \eqref{eq:ridgeprob_reg}, the approximation $\tilde{V}$ is given by 
\eqref{eq:approxValue} and the feedback law is recovered via \eqref{eq:feedhjb}.

\subsubsection{Partial use of second-order information}
\label{sec:partial_hessian}
Each sample contributes $n(n+1)/2$ distinct second-order equations, so the
second-order block dominates both the cost of data generation and the size of
the regression system as $n$ grows. We address this by subsampling the
second-order equations: given a ratio $\rho\in(0,1]$, we select a subset
$\mathcal{J}_\rho\subset\{1,\dots,N\}$ with $|\mathcal{J}_\rho|=\lfloor\rho N\rfloor$
and enforce Hessian equations only at those samples. The reduced stacked
operators $\bar{\mathbf{A}}_{\gamma_1,\gamma_2}^{(\rho)}$,
$\bar{\mathbf{V}}_{\gamma_1,\gamma_2}^{(\rho)}$ are assembled analogously to
\eqref{eq:weighted_system}, and problem \eqref{eq:weighted_leastsquares_compact}
is solved with total equation count
$M_\rho = N(1+n)+\lfloor\rho N\rfloor\,n(n+1)/2$.
Even moderate $\rho<1$ yields a substantial
reduction in both generation and regression costs.

\begin{remark}[Directional second-order projections]
An alternative strategy replaces the full Hessian by directional projections 
of the form $v^\top\nabla^2 V(\mathbf{x}^j)v$, with $v$ sampled from the unit 
sphere, following the Sobolev training framework of \cite{COJSP}. This retains 
global curvature information at reduced cost per sample. We do not pursue this 
variant here, focusing instead on subsampling at the level of sample points.
\end{remark}

\subsection{Performance and error measurement}
\label{sec:errors}

The generated dataset is split into a training set and an independent validation 
set. The training set is used to compute $\theta^*$, while the validation set 
is used exclusively to assess approximation accuracy. When testing partial 
second-order formulations, the subsampling 
strategy is applied only during training; all validation errors are computed 
using the complete first- and second-order information, ensuring a consistent 
comparison across models.

Let $\mathcal{I}_{\mathrm{val}}$ denote the validation index set and let 
$\tilde{V}$ be the fitted approximation obtained from 
\eqref{eq:weighted_leastsquares_compact} or  \eqref{eq:ridgeprob_reg}. Define the discrete $H^\ell$ seminorm
\begin{equation}\label{eq:disc_seminorm}
  \lvert W \rvert_{\ell,\mathrm{val}}^{2}
  \coloneqq
  \sum_{j\in\mathcal{I}_{\mathrm{val}}}
  \sum_{\lvert\alpha\rvert=\ell}
  \bigl\lvert D^{\boldsymbol{\alpha}}W(\mathbf{x}^{j})\bigr\rvert^{2},
  \qquad \ell = 0,1,2,
\end{equation}
where the inner sum runs over all multi-indices $\boldsymbol{\alpha}\in\mathbb{N}_0^n$ with 
$|\boldsymbol{\alpha}|=\ell$. The relative discrete Sobolev errors are
\begin{equation}\label{eq:errors}
  \mathrm{Err}_{H^{\ell}}(\tilde{V})
  \coloneqq
  \frac{\lvert \tilde{V} - V\rvert_{\ell,\mathrm{val}}}
       {\lvert V\rvert_{\ell,\mathrm{val}}},
  \qquad \ell=0,1,2,
\end{equation}
providing complementary information: $\ell=0$ quantifies global value 
approximation, $\ell=1$ reflects gradient accuracy and is directly related 
to feedback quality, and $\ell=2$ measures curvature consistency.

We complement the Sobolev errors with two metrics that 
assess approximation quality more directly in terms of the control problem.

The first is a PDE-based metric. For each validation point 
$\mathbf{x}^j\in\mathcal{I}_{\mathrm{val}}$, the pointwise HJB residual is
\begin{equation}\label{eq:hjb_residual}
r(\mathbf{x}^j) := \partial_t V(0,\mathbf{x}^j) + \ell(\mathbf{x}^j) 
+ \nabla\tilde{V}(\mathbf{x}^j)^\top\mathbf{f}(\mathbf{x}^j) 
- \frac{1}{4\beta}\|\mathbf{g}(\mathbf{x}^j)^\top\nabla\tilde{V}(\mathbf{x}^j)\|_2^2,
\end{equation}
where $\partial_t V(0,\mathbf{x}^j)$ is recovered by rearranging \eqref{HJBe} 
and substituting the exact values and gradient data available from $\mathcal{D}$. 
The advection and control terms in \eqref{eq:hjb_residual} are evaluated using 
the approximate gradient $\nabla\tilde{V}$, so the residual measures the extent 
to which $\tilde{V}$ fails to satisfy the HJB equation. Writing 
$\ell^j:=\ell(\mathbf{x}^j)$, $\mathbf{f}^j:=\mathbf{f}(\mathbf{x}^j)$, 
$\mathbf{g}^j:=\mathbf{g}(\mathbf{x}^j)$, $V_t^j:=\partial_tV(0,\mathbf{x}^j)$, 
and $\nabla V^j:=\nabla V(\mathbf{x}^j)$, the relative HJB residual is
\begin{equation}\label{eq:hjb_relative}
\mathrm{HJB}_{\mathrm{rel}}(\tilde{V}) := \frac{1}{|\mathcal{I}_{\mathrm{val}}|} 
\sum_{j \in \mathcal{I}_{\mathrm{val}}} 
\frac{|r(\mathbf{x}^j)|}{|V_t^j| + |\ell^j| 
+ |(\nabla V^j)^\top \mathbf{f}^j| 
+ \tfrac{1}{4\beta}\|(\mathbf{g}^j)^\top\nabla V^j\|_2^2}.
\end{equation}
The denominator normalises term by term using exact data from $\mathcal{D}$, 
yielding a scale-invariant measure that is not corrupted by approximation 
error in $\tilde{V}$.

The second metric quantifies the control error directly. The pointwise relative 
discrepancy between the approximate feedback \eqref{eq:feedhjb} applied to 
$\tilde{V}$ and the optimal control is
\begin{equation}\label{eq:ctrl_error}
e_{\mathrm{ctrl}}(\mathbf{x}^j) := 
\frac{\|\tilde{\mathbf{u}}(\mathbf{x}^j) - \mathbf{u}^*(\mathbf{x}^j)\|_2}
     {\|\mathbf{u}^*(\mathbf{x}^j)\|_2}.
\end{equation}
This quantity is uninformative near points where $\mathbf{u}^*(\mathbf{x}^j)\approx\mathbf{0}$; 
in the numerical experiments we therefore report spatial statistics of 
$e_{\mathrm{ctrl}}$ restricted to the subset of the validation set where 
$\|\mathbf{u}^*(\mathbf{x}^j)\|_2$ exceeds a prescribed threshold.

\FloatBarrier

\section{Numerical experiments}
\label{sec:numerics}

We first assess the proposed method on an analytical benchmark, where 
exact values, gradients, and Hessians are available in closed form, and then on three 
optimal control problems in dimensions $n=2$, $6$, and $19$, with training 
datasets generated by the PMP--Riccati pipeline. 
Each experiment compares the zeroth-, first-, second-, and partial second-order 
formulations, and states the sample budget in each subsection. All experiments were performed in MATLAB R2025a on an ASUS laptop with a
24-core Intel Core i9-14900HX CPU (2.2~GHz base) and 64~GB of RAM, running
Windows~11. The code reproducing the Hessian-augmented regression in all our tests is available at \href{https://github.com/mgomezaedo/HJB-Hessian-Learning}{https://github.com/mgomezaedo/HJB-Hessian-Learning}.

\subsection{Analytical benchmark}
\label{sec:analytical}
To isolate the approximation properties of the regression framework from the 
PMP--Riccati data-generation approach, we consider the function
\begin{equation}\label{eq:test_function_12d}
f_{AS}(\mathbf{y})
=
\prod_{j=1}^n 
\frac{n/4}{n/4 + (y_j - a_j)^2},
\qquad 
a_j = \frac{(-1)^j}{j+1},
\end{equation}
on $[-1,1]^{n}$, previously used in \cite{AS} to benchmark gradient-augmented 
polynomial approximation. We extend that comparison to second-order augmentation. We fix $n=12$.
The basis is $\mathfrak{I}_{\mathrm{HC}}(4)$, giving $q=3{,}482$ Legendre 
polynomials, with the domain extended to $[-1.15,1.15]^{12}$ following 
Remark~\ref{rem:rescaling}. Ridge regression is used with $\lambda=10^{-8}$ 
and $\gamma_1=\gamma_2=1$; training points are drawn from a scrambled Sobol 
sequence, and errors are averaged over $10$ independent training sets on a 
fixed validation set of $500$ points.

\begin{figure}[h!]
\centering
\captionsetup{aboveskip=-1pt,belowskip=-1pt}
\includegraphics[width=\textwidth]{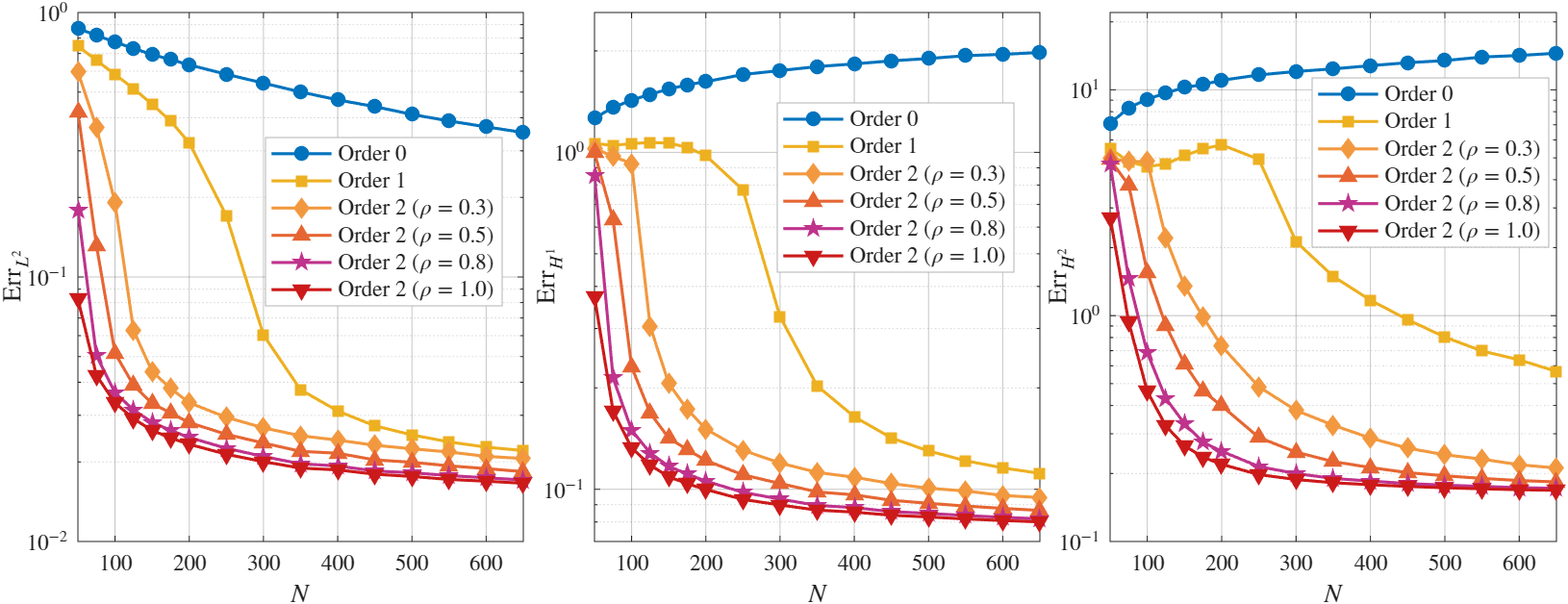} 
\caption{Relative $L^2$, $H^1$, and $H^2$ errors for the approximation of 
the 12-dimensional test \eqref{eq:test_function_12d} as a function 
of $N$, using $\mathfrak{I}_{\mathrm{HC}}(4)$ with $q = 3{,}482$. Zeroth-, 
first-, and second-order regressions are compared, with partial Hessian 
fractions $\rho \in \{0.3, 0.5, 0.8, 1.0\}$. Curves are averaged over $10$ 
independent training sets. Note that each panel reports a different error 
metric, identified by the subscript of $\mathrm{Err}$, and that the ordinate 
scales differ accordingly across panels.}
\label{fig:12d_ridge}
\end{figure}

Figure~\ref{fig:12d_ridge} illustrates the sample complexity gains from 
derivative augmentation. Each training point contributes $1$, $13$, or $91$ 
equations to the zeroth-, first-, and second-order systems respectively. 
Second-order augmentation achieves $L^2$ error below $7\%$ at $N=50$, a level 
first-order regression does not reach until $N\approx 350$; zeroth-order 
regression remains above $80\%$ throughout. Even $\rho=0.3$ consistently 
outperforms first-order regression across all metrics. The $H^2$ error reveals 
a persistent gap: without explicit Hessian data, curvature accuracy does not 
improve with $N$. We report that enriching the basis to $\mathfrak{I}_{\mathrm{HC}}(4)\cup 
\mathfrak{I}_{\mathrm{TD}}(6)$ ($q=18{,}948$) further illustrates this: at 
$N=700$, second-order augmentation reaches $0.74\%$ $L^2$ error, below the 
$1.5\%$ floor of the smaller basis, while first-order regression on the same 
basis stalls at $35\%$ with an underdetermined system.

We now turn to finite-horizon optimal control problems of the form 
\eqref{eq:optcost}, with datasets generated as described in 
Section~\ref{sec:augmented_data}.

\subsection{2-dimensional Van der Pol oscillator}
\label{sec:vdp_2d}

We consider the finite-horizon optimal control problem associated with the 
two-dimensional Van der Pol oscillator,
\begin{equation}\label{eq:vdp_cost}
\min_{u\in L^2([0,T];\mathbb{R})} 
\int_{0}^{T}y_1^2(t)+y_2^2(t)+\beta\,u^2(t)\,dt,
\end{equation}
subject to
\begin{equation}\label{eq:vdp_state}
\left\{
\begin{array}{l@{\qquad}l}
\partial_t y_1(t)=y_2(t), & y_1(0)=x_1,\\
\partial_t y_2(t)=-y_1(t)+y_2(t)\bigl(1-y_1^2(t)\bigr)+u(t), & y_2(0)=x_2.
\end{array}
\right.
\end{equation}
The associated adjoint variables $\mathbf{p}=(p_1,p_2)$ satisfy
\begin{equation}\label{eq:vdp_adjoint}
\left\{
\begin{array}{l@{\qquad}l}
\partial_t p_1(t)=\bigl(1+2y_1(t)y_2(t)\bigr)p_2(t)-2y_1(t), & p_1(T)=0,\\
\partial_t p_2(t)=-p_1(t)-\bigl(1-y_1^2(t)\bigr)p_2(t)-2y_2(t), & p_2(T)=0,
\end{array}
\right.
\end{equation}
and the optimal control is $u^*(t)=-\frac{1}{2\beta}\,p_2(t)$. The 
time-dependent coefficient matrices in the Riccati equation \eqref{eq:riccati} 
are given by
\begin{equation}\label{eq:vdp_riccati_coeffs}
\begin{array}{ll}
\mathcal{H}_{py}=
\begin{pmatrix}
0 & 1\\
-2y_1y_2-1 & 1-y_1^2
\end{pmatrix},
&
\mathcal{H}_{yp}=
\begin{pmatrix}
0 & -2y_1y_2-1\\
1 & 1-y_1^2
\end{pmatrix},
\\[10pt]
\mathcal{H}_{pp}=
\begin{pmatrix}
0 & 0\\
0 & -\frac{1}{2\beta}
\end{pmatrix},
&
\mathcal{H}_{yy}=
\begin{pmatrix}
2-2y_2p_2 & -2y_1p_2\\
-2y_1p_2 & 2
\end{pmatrix},
\end{array}
\end{equation}
with $\nabla^2 V(\mathbf{x})=R(0)$.

To enable a direct comparison with Test~1 in \cite{AKK}, we fix $T=3$ and 
$\beta=0.1$, and sample initial conditions uniformly from $\Omega=[-3,3]^2$. 
The approximation space is built from the hyperbolic-cross Legendre basis 
$\mathfrak{I}_{\mathrm{HC}}(4)$ ($q=52$), optionally enriched with total-degree 
terms $\mathfrak{I}_{\mathrm{TD}}(s_{\mathrm{td}})$ as in \eqref{eq:enriched_index}, 
yielding cardinalities up to $q=153$. Basis functions are rescaled to 
$[-3.45,3.45]^2$.

\begin{figure}[h!]
\centering
\begin{subfigure}[t]{0.32\textwidth}
    \centering
    \includegraphics[width=\linewidth]{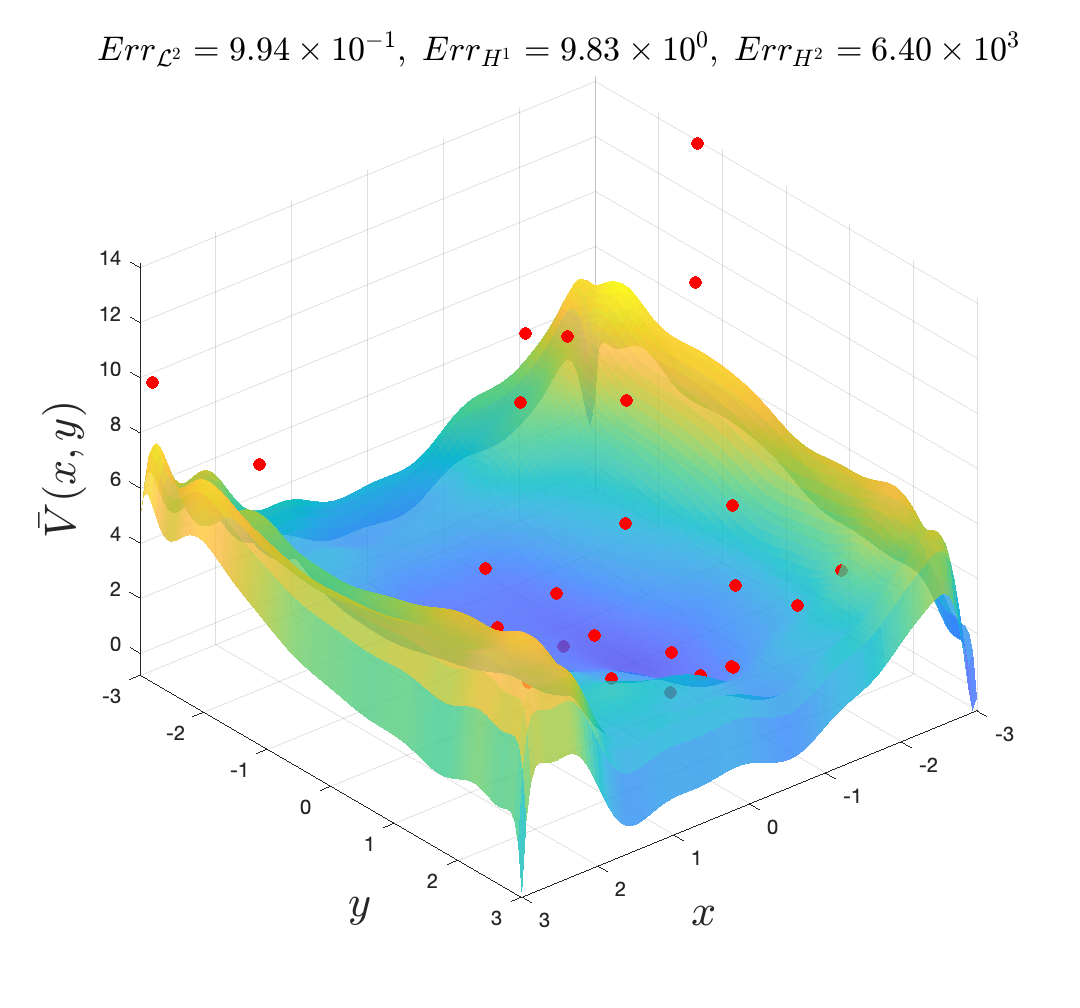}
    \subcaption*{(i) Zeroth-order}
\end{subfigure}\hfill
\begin{subfigure}[t]{0.32\textwidth}
    \centering
    \includegraphics[width=\linewidth]{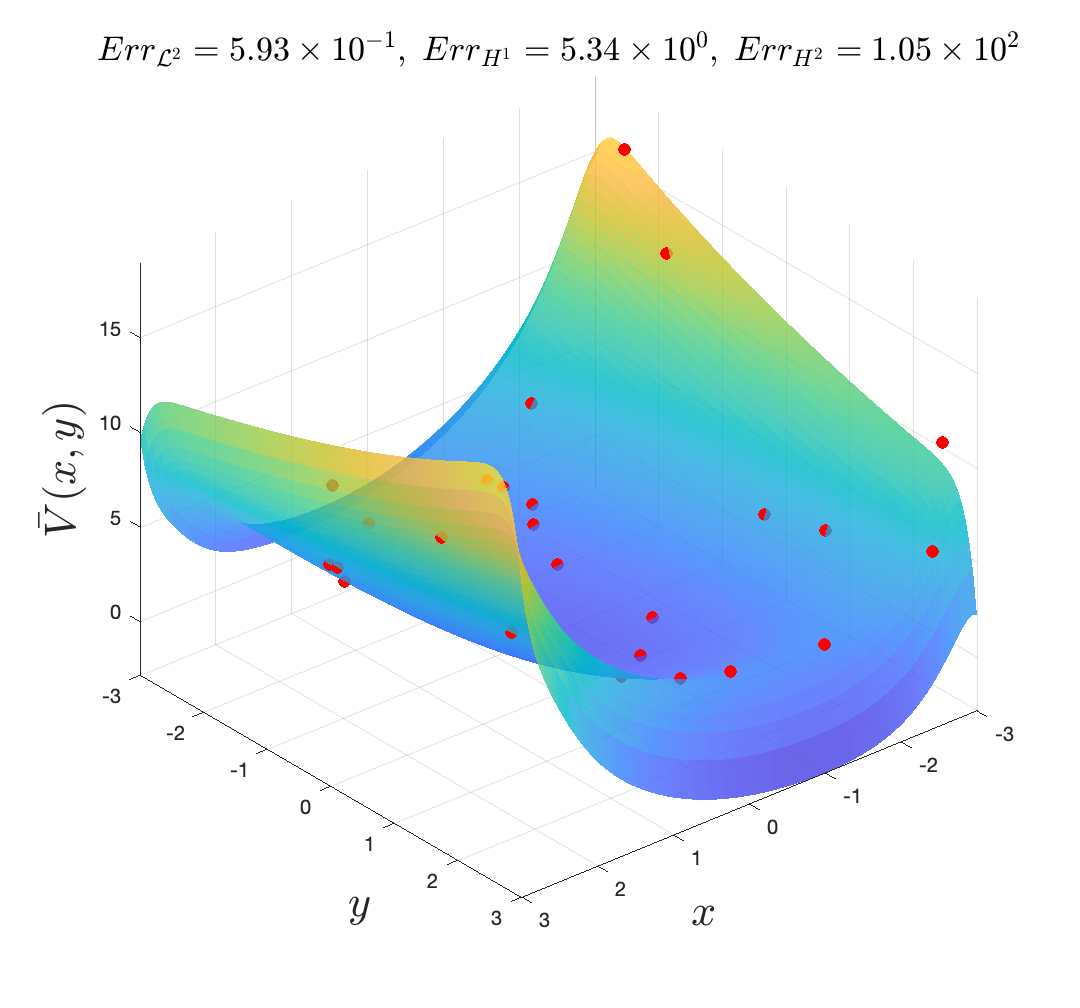}
    \subcaption*{(ii) First-order}
\end{subfigure}\hfill
\begin{subfigure}[t]{0.32\textwidth}
    \centering
    \includegraphics[width=\linewidth]{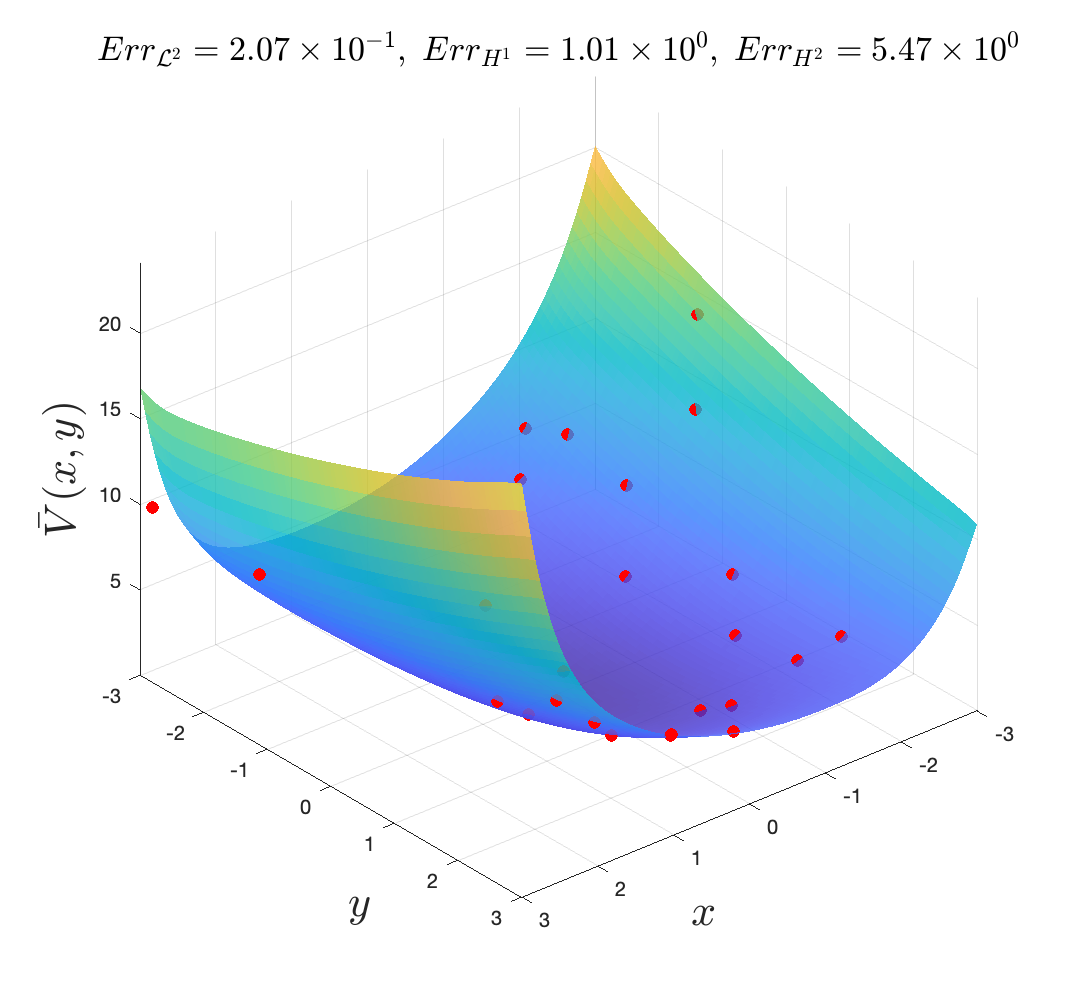}
    \subcaption*{(iii) Second-order}
\end{subfigure}
\caption{Polynomial approximation of the value function for the 2D Van der 
Pol problem with $N=25$ training samples (red dots) under zeroth-, first-, 
and second-order augmentation.}
\label{fig:vdp_surfaces}
\end{figure}

Figure~\ref{fig:vdp_surfaces} illustrates qualitatively the effect of 
derivative augmentation at $N=25$: zeroth-order regression produces a surface 
with pronounced oscillations, while second-order augmentation already yields 
a smooth, accurate reconstruction at the same sample count.

\begin{figure}[h!]
\centering
\begin{subfigure}[t]{\textwidth}
    \centering
    \includegraphics[width=\textwidth]{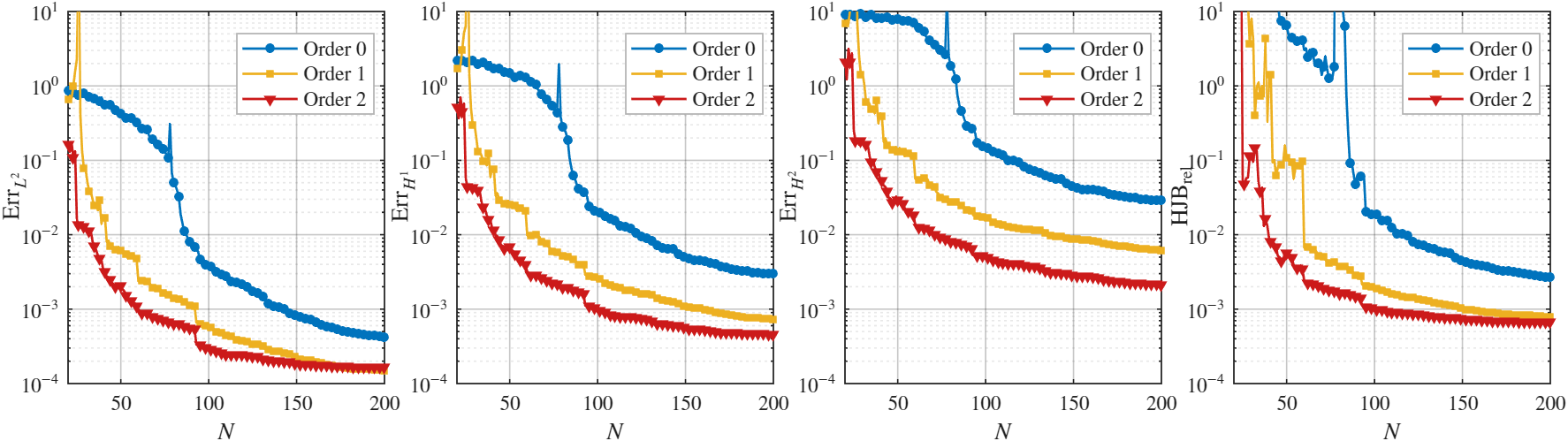}
    \caption{Least-squares regression}
    \label{fig:err_vs_N_ls}
\end{subfigure}
\\[0.3cm]
\begin{subfigure}[t]{\textwidth}
    \centering
    \includegraphics[width=\textwidth]{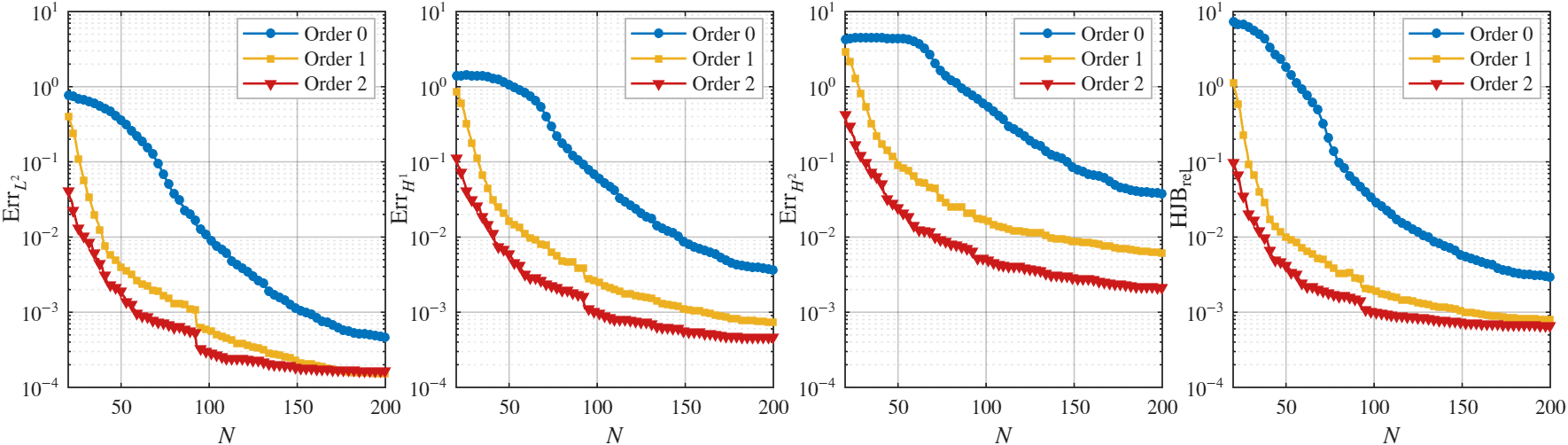}
    \caption{Ridge regression ($\lambda = 10^{-10}$)}
    \label{fig:err_vs_N_ridge}
\end{subfigure}
\caption{Relative validation errors as functions of $N$ for the enriched 
basis $\mathfrak{I}_{\mathrm{HC}}(4)\cup\mathfrak{I}_{\mathrm{TD}}(10)$ 
($q=78$), averaged over $100$ Monte Carlo realisations. 
Weights $\gamma_1=\gamma_2=1$. Note that each panel reports a different 
error metric, identified by the subscript of $\mathrm{Err}$, and that the 
ordinate scales differ accordingly across panels.}
\label{fig:err_vs_N}
\end{figure}

Figure~\ref{fig:err_vs_N} reports validation errors as functions of $N$ for 
the enriched basis $\mathfrak{I}_{\mathrm{HC}}(4)\cup\mathfrak{I}_{\mathrm{TD}}(10)$ 
($q=78$), averaged over $100$ Monte Carlo realisations. Panel~(a) uses standard 
least-squares; panel~(b) adds ridge regularisation with $\lambda=10^{-10}$, 
which removes instabilities in the underdetermined regime without affecting 
asymptotic accuracy. Second-order augmentation achieves the lowest errors 
across all four metrics at small to moderate $N$, with the $H^2$ gap 
persisting throughout regardless of sample size.

\begin{figure}[h!]
\centering
\includegraphics[width=\textwidth]{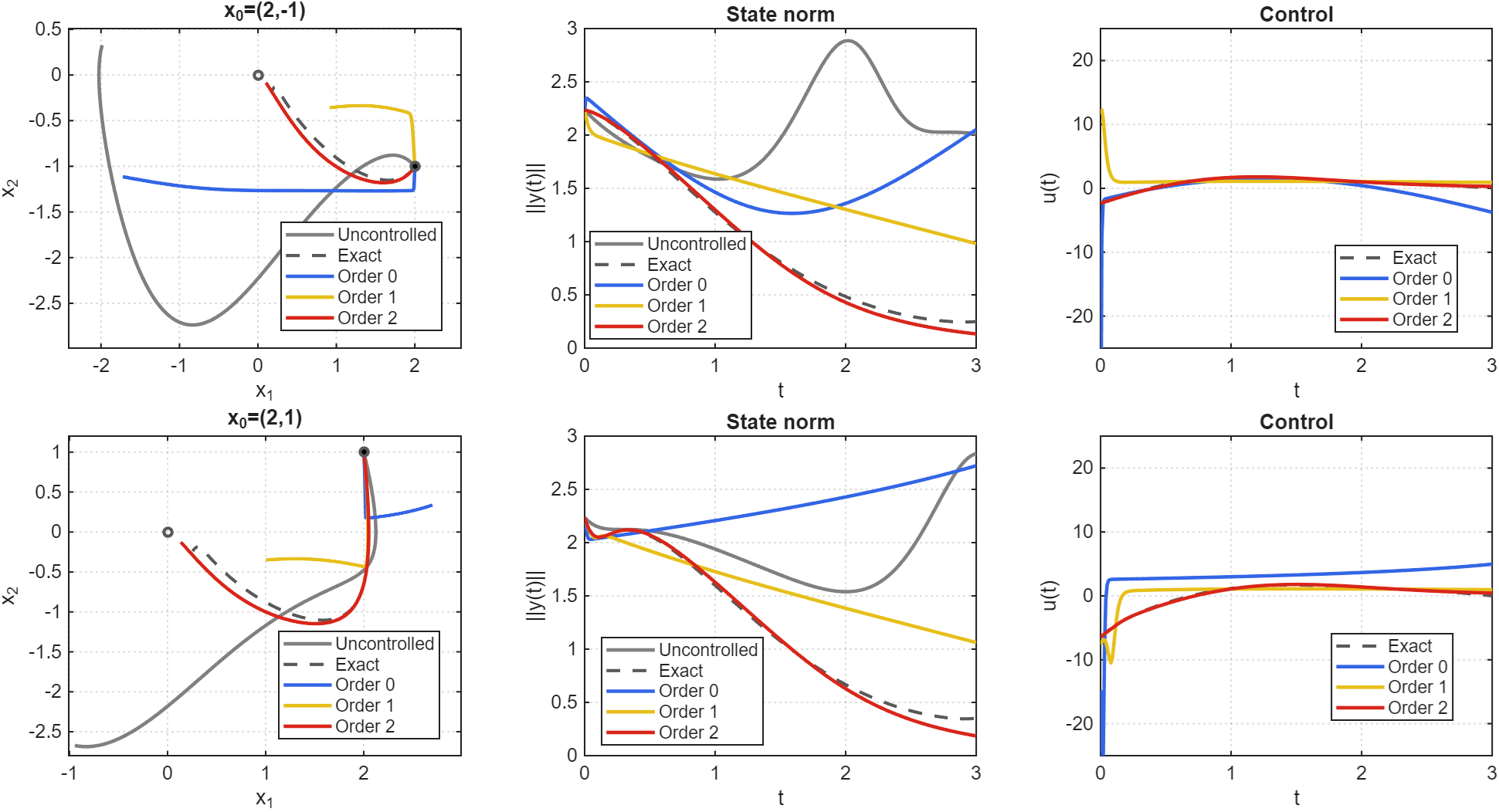}
\caption{Closed-loop trajectories for $N=15$ training samples. Each row 
corresponds to an initial condition $\mathbf{x}_0\in\{(2,-1),(2,1)\}$. 
Left: phase portrait. Centre: state norm $\|\mathbf{y}(t)\|_2$. Right: 
control signal $u(t)$, clipped to $[-25,25]$.}
\label{fig:traj_N15}
\end{figure}

Figure~\ref{fig:traj_N15} examines closed-loop performance at $N=15$, using 
the enriched basis $\mathfrak{I}_{\mathrm{HC}}(4)\cup\mathfrak{I}_{\mathrm{TD}}(8)$ 
($q=61$), ridge regularisation with $\lambda=10^{-10}$, and a single random 
training set. At this sample size, zeroth- and first-order feedbacks fail to 
stabilise the system for both initial conditions shown; second-order 
augmentation is the only formulation that produces trajectories closely 
tracking the PMP-optimal reference.

We next study how approximation quality depends jointly on basis dimension 
and derivative order. Fixing $N$, we vary $s_{\mathrm{td}}\in\{7,\ldots,16\}$ 
in $\mathfrak{I}_{\mathrm{HC}}(4)\cup\mathfrak{I}_{\mathrm{TD}}(s_{\mathrm{td}})$, 
giving $q\in[56,153]$, and include $\mathfrak{I}_{\mathrm{HC}}(4)$ ($q=52$) 
as a reference. Results are averaged over $100$ Monte Carlo realisations and 
include partial Hessian formulations with $\rho\in\{0.3,0.5\}$.

\begin{figure}[h!]
\centering
\includegraphics[width=\textwidth]{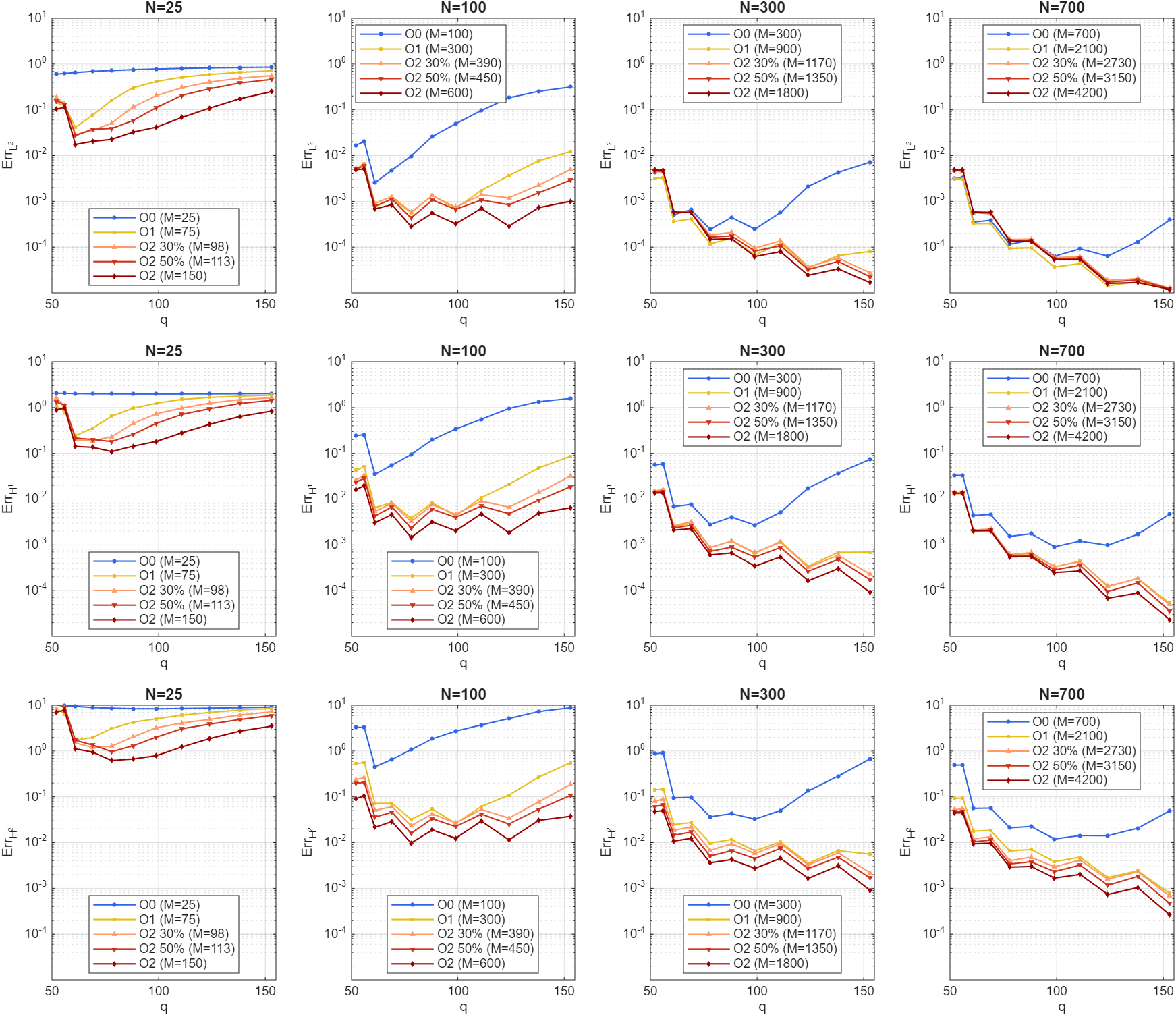}
\caption{Relative approximation errors as a function of basis cardinality $q$, 
for four values of $N$ and three error metrics. The total number of equations 
$M$ is indicated in each legend. Partial second-order formulations with 
$\rho\in\{0.3,0.5\}$ are included. Each point is the mean over $100$ Monte 
Carlo realisations; $\lambda=10^{-10}$, $\gamma_1=\gamma_2=1$. Note that 
each column reports a different error metric, identified by the subscript 
of $\mathrm{Err}$, and that the ordinate scales differ accordingly across 
panels.}
\label{fig:errvsq}
\end{figure}

Figure~\ref{fig:errvsq} shows the expected U-shaped behaviour in $q$ for 
each $(N,\text{Order})$ pair: errors decrease as the basis grows richer, 
then increase as the system loses overdetermination. The minimum shifts 
towards larger $q$ with derivative order, since higher-order augmentation 
maintains a well-conditioned system on more expressive bases, yielding lower 
approximation floors. Partial Hessian formulations interpolate smoothly 
between first- and second-order curves, with $\rho=0.3$ already providing 
a visible improvement over first-order regression across all metrics.

A practical design question is the minimum $N$ required to guarantee a 
prescribed tolerance $\tau$ on the relative control error. For each 
configuration $(N,\mathfrak{I},\text{Order})$ we perform $N_{\mathrm{MC}}=100$ 
independent Monte Carlo realisations. Each yields an approximate feedback law 
whose performance is summarised by the spatial $95$th percentile of the 
pointwise control error \eqref{eq:ctrl_error} over the active domain 
$\mathcal{I}_{\mathrm{act}}=\{j\in\mathcal{I}_{\mathrm{val}}:
\|\mathbf{u}^*(\mathbf{x}^j)\|_2\geq 0.1\}$ ($99.4\%$ of validation points):
\begin{equation}\label{eq:p95_per_run}
P^{(i)} := \mathrm{p}_{95}\!\Bigl(
\bigl\{e_{\mathrm{ctrl}}^{(i)}(\mathbf{x}^j)\bigr\}_{j\in\mathcal{I}_{\mathrm{act}}}
\Bigr),
\qquad i=1,\ldots,N_{\mathrm{MC}}.
\end{equation}
The design criterion requires both the mean and median of 
$\{P^{(i)}\}_{i=1}^{N_{\mathrm{MC}}}$ to lie below $\tau$:
\begin{equation}\label{eq:design_criterion}
\frac{1}{N_{\mathrm{MC}}}\sum_{i=1}^{N_{\mathrm{MC}}}P^{(i)}\leq\tau
\qquad\text{and}\qquad
\mathrm{median}\!\left(P^{(1)},\ldots,P^{(N_{\mathrm{MC}})}\right)\leq\tau,
\end{equation}
ensuring that a representative realisation meets the tolerance while 
preventing outliers from being concealed in the average. For each derivative 
order we search over bases $\mathfrak{I}_{\mathrm{HC}}(4)\cup
\mathfrak{I}_{\mathrm{TD}}(s_{\mathrm{td}})$ with $s_{\mathrm{td}}\in\{7,\ldots,16\}$ 
and record the minimum $N$ satisfying \eqref{eq:design_criterion}. 
Table~\ref{tab:min_N_comparison} reports the outcome: second-order 
augmentation consistently requires the fewest samples, with the ratio 
$N_{\mathrm{O0}}/N_{\mathrm{O2}}$ growing from $3.5$ at $\tau=10\%$ to 
$11.7$ at $\tau=0.1\%$.

\begin{table}[h!]
\centering
\small
\setlength{\tabcolsep}{4pt}
\resizebox{\textwidth}{!}{%
\begin{tabular}{c c c c c c c c c}
\toprule
$\tau$ & $N_{\mathrm{O2}}$ & $N_{\mathrm{O1}}$ &
$N_{\mathrm{O0}}$ & $\frac{N_{\mathrm{O1}}}{N_{\mathrm{O2}}}$
& $\frac{N_{\mathrm{O0}}}{N_{\mathrm{O2}}}$ &
Basis O2 & Basis O1 & Basis O0 \\
\midrule
$10\%$  & $\mathbf{29}$  & $48$   & $102$  & $1.7$ & $3.5$
& $\mathfrak{I}_{\mathrm{HC}}(4)\!\cup\!\mathfrak{I}_{\mathrm{TD}}(8)$
& $\mathfrak{I}_{\mathrm{HC}}(4)\!\cup\!\mathfrak{I}_{\mathrm{TD}}(10)$
& $\mathfrak{I}_{\mathrm{HC}}(4)\!\cup\!\mathfrak{I}_{\mathrm{TD}}(8)$  \\
$5\%$   & $\mathbf{38}$  & $61$   & $138$  & $1.6$ & $3.6$
& $\mathfrak{I}_{\mathrm{HC}}(4)\!\cup\!\mathfrak{I}_{\mathrm{TD}}(8)$
& $\mathfrak{I}_{\mathrm{HC}}(4)\!\cup\!\mathfrak{I}_{\mathrm{TD}}(10)$
& $\mathfrak{I}_{\mathrm{HC}}(4)\!\cup\!\mathfrak{I}_{\mathrm{TD}}(8)$  \\
$1\%$   & $\mathbf{80}$  & $109$  & $277$  & $1.4$ & $3.5$
& $\mathfrak{I}_{\mathrm{HC}}(4)\!\cup\!\mathfrak{I}_{\mathrm{TD}}(14)$
& $\mathfrak{I}_{\mathrm{HC}}(4)\!\cup\!\mathfrak{I}_{\mathrm{TD}}(12)$
& $\mathfrak{I}_{\mathrm{HC}}(4)\!\cup\!\mathfrak{I}_{\mathrm{TD}}(12)$ \\
$0.5\%$ & $\mathbf{87}$  & $134$  & $372$  & $1.5$ & $4.3$
& $\mathfrak{I}_{\mathrm{HC}}(4)\!\cup\!\mathfrak{I}_{\mathrm{TD}}(14)$
& $\mathfrak{I}_{\mathrm{HC}}(4)\!\cup\!\mathfrak{I}_{\mathrm{TD}}(12)$
& $\mathfrak{I}_{\mathrm{HC}}(4)\!\cup\!\mathfrak{I}_{\mathrm{TD}}(12)$ \\
$0.1\%$ & $\mathbf{144}$ & $208$  & $1689$ & $1.4$ & $11.7$
& $\mathfrak{I}_{\mathrm{HC}}(4)\!\cup\!\mathfrak{I}_{\mathrm{TD}}(14)$
& $\mathfrak{I}_{\mathrm{HC}}(4)\!\cup\!\mathfrak{I}_{\mathrm{TD}}(14)$
& $\mathfrak{I}_{\mathrm{HC}}(4)\!\cup\!\mathfrak{I}_{\mathrm{TD}}(14)$ \\
\bottomrule
\end{tabular}}
\caption{Minimum number of training samples required by each derivative order 
to satisfy the control tolerance $\tau$ under criterion \eqref{eq:design_criterion}, 
with basis selected independently per order to minimise $N$. The ratios 
$N_{\mathrm{O1}}/N_{\mathrm{O2}}$ and $N_{\mathrm{O0}}/N_{\mathrm{O2}}$ 
quantify the data savings provided by second-order augmentation.}
\label{tab:min_N_comparison}
\end{table}

Table~\ref{tab:tolerance_mean} fixes $N$ at the second-order minimum and 
reports the best performance attainable by lower-order formulations at that 
same budget, with basis selected independently per order. For zeroth-order, 
bases with $q\leq N$ are preferred; when none exists in the family, the 
smallest available $q$ is used and the entry is marked with an asterisk. 
Second-order augmentation is the only formulation satisfying 
\eqref{eq:design_criterion} at every tolerance; the ratios $M_k/q$ explain 
the mechanism: only second-order augmentation maintains a well-determined 
system on a sufficiently expressive basis. In a typical realisation the 
spatial $\mathrm{p}_{95}$ is $2$--$4$ times smaller than the mean reported 
here, so the design criterion is conservative.

\begin{table}[h!]
\centering
\footnotesize
\setlength{\tabcolsep}{4pt}
\resizebox{\textwidth}{!}{%
\begin{tabular}{c c l c c c c c c c c c}
\toprule
$\tau$ & $N$ & Order & Basis & $q$ & $M$ &
$\frac{M_0}{q}$ & $\frac{M_1}{q}$ & $\frac{M_2}{q}$ &
$\frac{M}{q}$ & $\mathrm{p}_{95}$ & $\max$ \\
\midrule
\multirow{3}{*}{$10\%$}
& \multirow{3}{*}{$29$}
& Order~2 & $\mathfrak{I}_{\mathrm{HC}}(4)\!\cup\!\mathfrak{I}_{\mathrm{TD}}(8)$ & $61$ & $174$ &
  $0.5$ & $1.0$ & $1.4$ & $2.9$ &
  $9.7 \times 10^{-2}$ & $1.2 \times 10^{1}$ \\
& & Order~1 & $\mathfrak{I}_{\mathrm{HC}}(4)\!\cup\!\mathfrak{I}_{\mathrm{TD}}(8)$ & $61$ & $87$ &
  $0.5$ & $1.0$ & --- & $1.4$ &
  $3.9 \times 10^{-1}$ & $4.0 \times 10^{1}$ \\
& & Order~0$^*$ & $\mathfrak{I}_{\mathrm{HC}}(4)$ & $52$ & $29$ &
  $0.6$ & --- & --- & $0.6$ &
  $1.7 \times 10^{1}$ & $2.6 \times 10^{2}$ \\
\midrule
\multirow{3}{*}{$5\%$}
& \multirow{3}{*}{$38$}
& Order~2 & $\mathfrak{I}_{\mathrm{HC}}(4)\!\cup\!\mathfrak{I}_{\mathrm{TD}}(8)$ & $61$ & $228$ &
  $0.6$ & $1.2$ & $1.9$ & $3.7$ &
  $4.8 \times 10^{-2}$ & $3.8 \times 10^{0}$ \\
& & Order~1 & $\mathfrak{I}_{\mathrm{HC}}(4)\!\cup\!\mathfrak{I}_{\mathrm{TD}}(8)$ & $61$ & $114$ &
  $0.6$ & $1.2$ & --- & $1.9$ &
  $2.0 \times 10^{-1}$ & $2.2 \times 10^{1}$ \\
& & Order~0$^*$ & $\mathfrak{I}_{\mathrm{HC}}(4)$ & $52$ & $38$ &
  $0.7$ & --- & --- & $0.7$ &
  $8.6 \times 10^{0}$ & $2.2 \times 10^{2}$ \\
\midrule
\multirow{3}{*}{$1\%$}
& \multirow{3}{*}{$80$}
& Order~2 & $\mathfrak{I}_{\mathrm{HC}}(4)\!\cup\!\mathfrak{I}_{\mathrm{TD}}(14)$ & $124$ & $480$ &
  $0.6$ & $1.3$ & $1.9$ & $3.9$ &
  $8.1 \times 10^{-3}$ & $5.3 \times 10^{-1}$ \\
& & Order~1 & $\mathfrak{I}_{\mathrm{HC}}(4)\!\cup\!\mathfrak{I}_{\mathrm{TD}}(10)$ & $78$ & $240$ &
  $1.0$ & $2.1$ & --- & $3.1$ &
  $2.2 \times 10^{-2}$ & $9.3 \times 10^{-1}$ \\
& & Order~0 & $\mathfrak{I}_{\mathrm{HC}}(4)\!\cup\!\mathfrak{I}_{\mathrm{TD}}(8)$ & $61$ & $80$ &
  $1.3$ & --- & --- & $1.3$ &
  $2.4 \times 10^{-1}$ & $2.1 \times 10^{1}$ \\
\midrule
\multirow{3}{*}{$0.5\%$}
& \multirow{3}{*}{$87$}
& Order~2 & $\mathfrak{I}_{\mathrm{HC}}(4)\!\cup\!\mathfrak{I}_{\mathrm{TD}}(14)$ & $124$ & $522$ &
  $0.7$ & $1.4$ & $2.1$ & $4.2$ &
  $4.4 \times 10^{-3}$ & $2.2 \times 10^{-1}$ \\
& & Order~1 & $\mathfrak{I}_{\mathrm{HC}}(4)\!\cup\!\mathfrak{I}_{\mathrm{TD}}(10)$ & $78$ & $261$ &
  $1.1$ & $2.2$ & --- & $3.3$ &
  $2.4 \times 10^{-2}$ & $1.4 \times 10^{0}$ \\
& & Order~0 & $\mathfrak{I}_{\mathrm{HC}}(4)\!\cup\!\mathfrak{I}_{\mathrm{TD}}(8)$ & $61$ & $87$ &
  $1.4$ & --- & --- & $1.4$ &
  $1.7 \times 10^{-1}$ & $1.7 \times 10^{1}$ \\
\midrule
\multirow{3}{*}{$0.1\%$}
& \multirow{3}{*}{$144$}
& Order~2 & $\mathfrak{I}_{\mathrm{HC}}(4)\!\cup\!\mathfrak{I}_{\mathrm{TD}}(14)$ & $124$ & $864$ &
  $1.2$ & $2.3$ & $3.5$ & $7.0$ &
  $9.9 \times 10^{-4}$ & $3.0 \times 10^{-2}$ \\
& & Order~1 & $\mathfrak{I}_{\mathrm{HC}}(4)\!\cup\!\mathfrak{I}_{\mathrm{TD}}(12)$ & $99$ & $432$ &
  $1.5$ & $2.9$ & --- & $4.4$ &
  $4.2 \times 10^{-3}$ & $1.4 \times 10^{-1}$ \\
& & Order~0 & $\mathfrak{I}_{\mathrm{HC}}(4)\!\cup\!\mathfrak{I}_{\mathrm{TD}}(8)$ & $61$ & $144$ &
  $2.4$ & --- & --- & $2.4$ &
  $4.6 \times 10^{-2}$ & $5.4 \times 10^{0}$ \\
\bottomrule
\end{tabular}}
\caption{Best performance attainable by each derivative order at the minimum 
$N$ required by second-order augmentation to satisfy \eqref{eq:design_criterion}. 
Basis selected independently per order to minimise $\mathrm{p}_{95}$ of the 
mean control error. Asterisk: $q>N$. The ratios $M_k/q$ give equations per 
unknown at each derivative order; $M/q$ is the overall overdetermination ratio.}
\label{tab:tolerance_mean}
\end{table}

\begin{remark}[Computational cost]
The offline stage, which includes solving the BVP \eqref{eq:vdp_state}--\eqref{eq:vdp_adjoint} 
together with the Riccati equation \eqref{eq:riccati} for each sampled initial 
condition, takes a median of $93\,\mathrm{ms}$ per sample on the hardware 
described above, so a full training set is generated in seconds to a few minutes 
for the regimes considered here. Online evaluation of $\tilde{V}$, 
$\nabla\tilde{V}$, and $\nabla^2\tilde{V}$ requires only basis-function 
evaluations and dot products of length $q$, taking $0.5$--$1\,\mathrm{ms}$ 
across the basis sizes considered.
\end{remark}

\subsection{6-dimensional rigid body model of a satellite}
\label{sec:sat_6d}

We consider a six-state rigid body model of a satellite studied in \cite{NGK}. Here, $\mathbf{y}=(y_1,\dots,y_6)^\top\in\mathbb{R}^6$ 
is decomposed into attitude angles $\mathbf{y}_r=(y_1,y_2,y_3)^\top$ and 
angular velocity $\mathbf{y}_\omega=(y_4,y_5,y_6)^\top$, with 
torque input $\mathbf{u}\in\mathbb{R}^3$. The optimal control problem is
\begin{equation}\label{eq:sat_cost}
\min_{\mathbf{u}\in L^2([0,T];\mathbb{R}^3)}
\int_{0}^{T}\alpha_1\|\mathbf{y}_r\|_2^2+\alpha_2\|\mathbf{y}_\omega\|_2^2
+\beta\|\mathbf{u}\|_2^2\,dt+\alpha_3\|\mathbf{y}(T)\|_2^2,
\end{equation}
subject to
\begin{equation}\label{eq:sat_state}
\begin{cases}
\partial_t\mathbf{y}_r = \mathbf{E}(\mathbf{y}_r)\mathbf{y}_\omega,\\
\partial_t\mathbf{y}_\omega = \mathbf{J}\mathbf{S}(\mathbf{y}_\omega)
\mathbf{R}(\mathbf{y}_r)\mathbf{h}+\mathbf{J}\mathbf{B}\mathbf{u},
\end{cases}
\qquad\mathbf{y}(0)=\mathbf{x},
\end{equation}
where $\alpha_1,\alpha_2,\alpha_3,\beta>0$ and
{\small
\setlength{\arraycolsep}{3pt}
\[
\mathbf{E}(\mathbf{y}_r):=\begin{pmatrix}
1 & \sin y_1\tan y_2 & \cos y_1\tan y_2 \\
0 & \cos y_1 & -\sin y_1 \\
0 & \sin y_1/\cos y_2 & \cos y_1/\cos y_2
\end{pmatrix},
\qquad
\mathbf{S}(\mathbf{y}_\omega):=\begin{pmatrix}
0 & y_6 & -y_5 \\
-y_6 & 0 & y_4 \\
y_5 & -y_4 & 0
\end{pmatrix},
\]
\[
\mathbf{R}(\mathbf{y}_r):=\begin{pmatrix}
\cos y_2\cos y_3 & \cos y_2\sin y_3 & -\sin y_2 \\
\sin y_1\sin y_2\cos y_3-\cos y_1\sin y_3 & \sin y_1\sin y_2\sin y_3
+\cos y_1\cos y_3 & \cos y_2\sin y_1 \\
\cos y_1\sin y_2\cos y_3+\sin y_1\sin y_3 & \cos y_1\sin y_2\sin y_3
-\sin y_1\cos y_3 & \cos y_2\cos y_1
\end{pmatrix}\!\!,
\]
}
{\small
\[
\mathbf{J}:=\begin{pmatrix}
\tfrac{1}{2} & 0 & 0 \\
0 & \tfrac{1}{3} & 0 \\
0 & 0 & \tfrac{1}{4}
\end{pmatrix},
\qquad
\mathbf{h}:=\begin{pmatrix}1\\1\\1\end{pmatrix},
\qquad
\mathbf{B}:=\begin{pmatrix}
1 & \tfrac{1}{20} & \tfrac{1}{10} \\
\tfrac{1}{15} & 1 & \tfrac{1}{10} \\
\tfrac{1}{10} & \tfrac{1}{15} & 1
\end{pmatrix}.
\]
}
The terminal cost yields the boundary condition $\mathbf{p}(T)=2\alpha_3\,\mathbf{y}(T)$, 
and the stationarity condition gives 
$\mathbf{u}^*(t)=-\frac{1}{2\beta}\mathbf{B}^\top\mathbf{J}^\top
(p_4(t),p_5(t),p_6(t))^\top$; the full adjoint system is omitted for brevity.

We follow the parameter setting of \cite{NGK}: $\alpha_1=\tfrac{1}{2}$, 
$\alpha_2=5$, $\alpha_3=\tfrac{1}{2}$, $\beta=\tfrac{1}{4}$, $T=20$. 
Initial conditions are sampled from 
$\Omega=[-\pi/3,\pi/3]^3\times[-\pi/4,\pi/4]^3$ using a scrambled Sobol 
sequence. The TPBVP \eqref{eq:tpbvp} is solved with \texttt{bvp5c} 
($\mathrm{AbsTol}=10^{-8}$, $\mathrm{RelTol}=10^{-7}$). Bases are of the 
enriched form \eqref{eq:enriched_index} with $\mu=0.15$, $\lambda=10^{-8}$, 
and $\gamma_1=\gamma_2=1$. A fixed validation set of $500$ points is used; 
errors are averaged over $50$ Monte Carlo realisations.

\begin{figure}[h!]
\centering
\includegraphics[width=\textwidth]{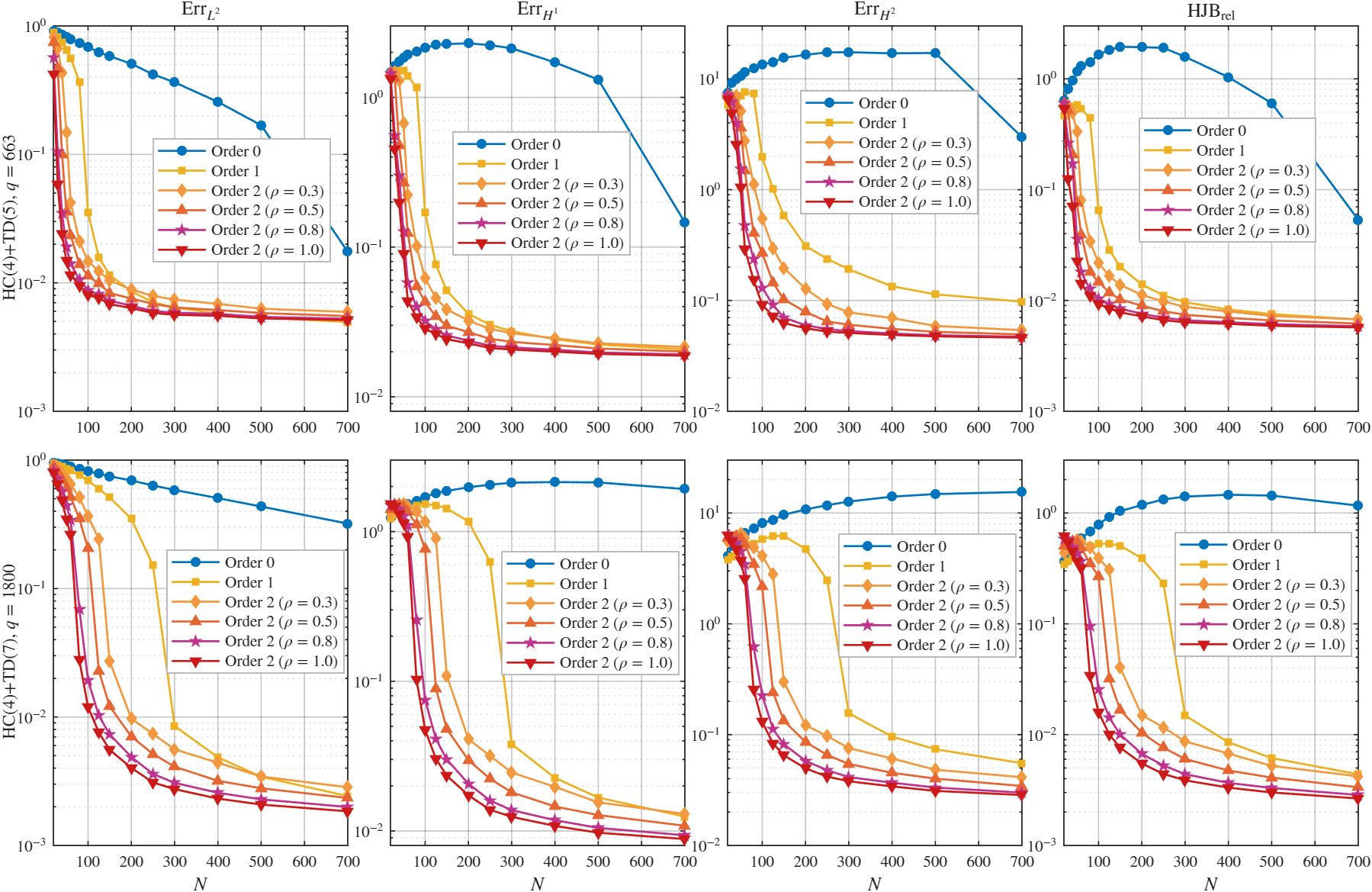}
\caption{Relative $L^2$, $H^1$, $H^2$ errors and HJB residual for the 6D 
satellite problem as a function of $N$, for bases 
$\mathfrak{I}_{\mathrm{HC}}(4)\cup\mathfrak{I}_{\mathrm{TD}}(5)$ ($q=663$, 
top row) and $\mathfrak{I}_{\mathrm{HC}}(4)\cup\mathfrak{I}_{\mathrm{TD}}(7)$ 
($q=1{,}800$, bottom row). Hessian fractions $\rho\in\{0.3,0.5,0.8,1.0\}$ 
included for second-order regression. Curves averaged over $50$ Monte Carlo 
realisations. Note that each panel reports a different error metric, 
identified by the subscript of $\mathrm{Err}$ (or the HJB residual), and 
that the ordinate scales differ accordingly across panels.}
\label{fig:sat6d_errors}
\end{figure}

Figure~\ref{fig:sat6d_errors} illustrates two qualitatively distinct regimes. 
On the moderate basis ($q=663$, top row), all formulations are feasible from 
small $N$; second-order augmentation reaches $Err_{L^2}\approx 10^{-2}$ near 
$N=60$, while first-order regression requires $N\approx 100$ to match the same 
accuracy. On the enriched basis ($q=1{,}800$, bottom row), the first-order 
system remains underdetermined ($M/q<1$) for $N\lesssim 250$, keeping 
$Err_{L^2}$ above $10^{-1}$ throughout that range; second-order augmentation, 
contributing $28N$ equations, reaches a well-determined regime at $N=100$ and 
achieves $Err_{L^2}=3.9\times 10^{-3}$ at $N=200$. As in previous experiments, 
partial Hessian formulations interpolate smoothly between the first- and 
second-order curves, with $\rho=0.3$ already providing a visible improvement.

To examine how far the basis dimension can be pushed with full second-order 
regression, we sweep $\mathfrak{I}_{\mathrm{HC}}(4)\cup\mathfrak{I}_{\mathrm{TD}}(s)$ 
for $s\in\{5,\ldots,11\}$ ($q$ from $663$ to $12{,}406$), restricting to 
$N\geq q/28$ to maintain a determined system, and select for each $N$ the 
basis minimising $Err_{L^2}$ averaged over $30$ Monte Carlo realisations.

\begin{table}[h!]
\centering
\small
\setlength{\tabcolsep}{4pt}
\resizebox{\textwidth}{!}{%
\begin{tabular}{rcccrcccc}
\toprule
$N$ & Best basis & $q$ & $M/q$ & & $Err_{L^2}$ & $Err_{H^1}$ & 
$Err_{H^2}$ & $\mathrm{HJB}_{\mathrm{rel}}$ \\
\midrule
$50$   & $\mathfrak{I}_{\mathrm{HC}}(4)\!\cup\!\mathfrak{I}_{\mathrm{TD}}(5)$  
& $663$    & $2.1$ & & $1.64\times10^{-2}$ & $1.06\times10^{-1}$ 
& $1.24$ & $4.19\times10^{-2}$ \\
$100$  & $\mathfrak{I}_{\mathrm{HC}}(4)\!\cup\!\mathfrak{I}_{\mathrm{TD}}(6)$  
& $1{,}044$  & $2.7$ & & $7.40\times10^{-3}$ & $2.79\times10^{-2}$ 
& $8.70\times10^{-2}$ & $9.06\times10^{-3}$ \\
$200$  & $\mathfrak{I}_{\mathrm{HC}}(4)\!\cup\!\mathfrak{I}_{\mathrm{TD}}(7)$  
& $1{,}800$  & $3.1$ & & $3.90\times10^{-3}$ & $1.67\times10^{-2}$ 
& $4.82\times10^{-2}$ & $5.34\times10^{-3}$ \\
$300$  & $\mathfrak{I}_{\mathrm{HC}}(4)\!\cup\!\mathfrak{I}_{\mathrm{TD}}(7)$  
& $1{,}800$  & $4.7$ & & $2.70\times10^{-3}$ & $1.24\times10^{-2}$ 
& $3.81\times10^{-2}$ & $3.80\times10^{-3}$ \\
$400$  & $\mathfrak{I}_{\mathrm{HC}}(4)\!\cup\!\mathfrak{I}_{\mathrm{TD}}(8)$  
& $3{,}051$  & $3.7$ & & $2.15\times10^{-3}$ & $1.06\times10^{-2}$ 
& $3.49\times10^{-2}$ & $3.10\times10^{-3}$ \\
$500$  & $\mathfrak{I}_{\mathrm{HC}}(4)\!\cup\!\mathfrak{I}_{\mathrm{TD}}(9)$  
& $5{,}047$  & $2.8$ & & $1.74\times10^{-3}$ & $1.00\times10^{-2}$ 
& $3.76\times10^{-2}$ & $2.96\times10^{-3}$ \\
$700$  & $\mathfrak{I}_{\mathrm{HC}}(4)\!\cup\!\mathfrak{I}_{\mathrm{TD}}(9)$  
& $5{,}047$  & $3.9$ & & $1.32\times10^{-3}$ & $7.22\times10^{-3}$ 
& $2.79\times10^{-2}$ & $2.03\times10^{-3}$ \\
$1000$ & $\mathfrak{I}_{\mathrm{HC}}(4)\!\cup\!\mathfrak{I}_{\mathrm{TD}}(10)$ 
& $8{,}044$  & $3.5$ & & $9.10\times10^{-4}$ & $5.61\times10^{-3}$ 
& $2.47\times10^{-2}$ & $1.54\times10^{-3}$ \\
$1500$ & $\mathfrak{I}_{\mathrm{HC}}(4)\!\cup\!\mathfrak{I}_{\mathrm{TD}}(11)$ 
& $12{,}406$ & $3.4$ & & $6.31\times10^{-4}$ & $4.35\times10^{-3}$ 
& $2.17\times10^{-2}$ & $1.14\times10^{-3}$ \\
$1800$ & $\mathfrak{I}_{\mathrm{HC}}(4)\!\cup\!\mathfrak{I}_{\mathrm{TD}}(11)$ 
& $12{,}406$ & $4.1$ & & $\mathbf{5.07\times10^{-4}}$ & 
$\mathbf{3.49\times10^{-3}}$ & $\mathbf{1.81\times10^{-2}}$ & 
$\mathbf{9.43\times10^{-4}}$ \\
\bottomrule
\end{tabular}}
\caption{Best basis of the form $\mathfrak{I}_{\mathrm{HC}}(4)\cup
\mathfrak{I}_{\mathrm{TD}}(s)$ for each $N$, selected to minimise $Err_{L^2}$ 
over $s\in\{5,\ldots,11\}$. The ratio $M/q=28N/q$ measures the degree of 
overdetermination. All errors averaged over $30$ Monte Carlo realisations.}
\label{tab:6d_best_basis}
\end{table}

The optimal overdetermination ratio $M/q$ lies consistently in the range 
$3$--$5$, confirming that second-order augmentation supports increasingly 
expressive bases as $N$ grows. At the largest configuration ($N=1{,}800$, 
$q=12{,}406$), lower-order regressions are not viable: the first-order system 
has $M/q\approx 1.02$ and the zeroth-order system $M/q\approx 0.15$, both 
effectively rank-deficient at this basis size.

\begin{remark}
The $Err_{H^2}$ values in Table~\ref{tab:6d_best_basis} remain an order of 
magnitude above the other metrics. This is expected: the Hessian block acts 
primarily as a regulariser that transfers second-order information to the 
coefficient vector $\theta$, improving $\tilde{V}$, $\nabla\tilde{V}$, and 
the feedback law \eqref{eq:feedhjb}; driving $Err_{H^2}$ to zero would 
overfit curvature at the expense of the metrics most relevant for control.
\end{remark}

\begin{figure}[h!]
\centering
\includegraphics[width=\linewidth]{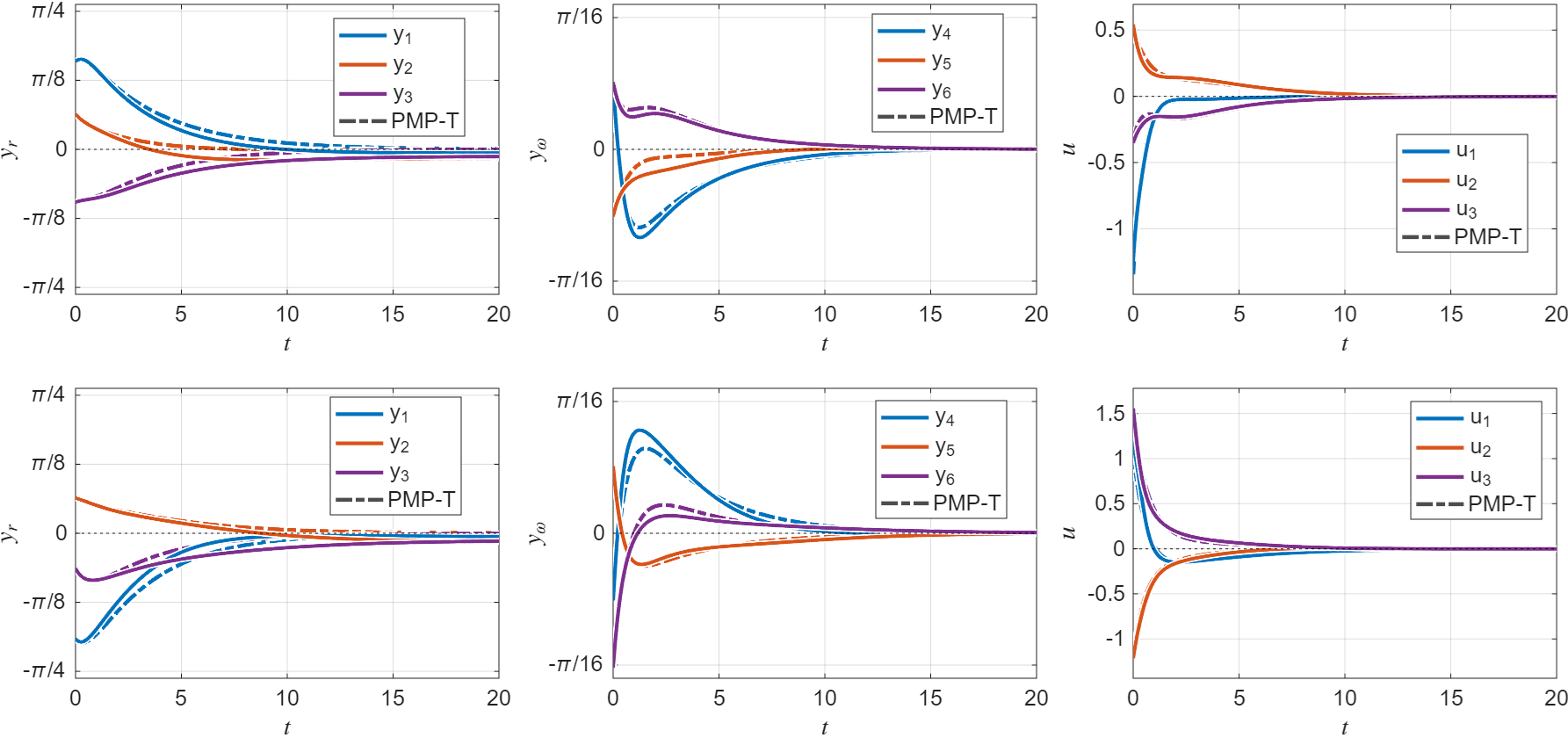}
\caption{Closed-loop trajectories for the 6D satellite using partial Hessian
regression ($\rho=0.5$, $N=70$, $\mathfrak{I}_{\mathrm{HC}}(4)$, $q=537$),
from $\mathrm{IC}_1=(0.5,0.2,-0.3,0.1,-0.1,0.1)$ (top) and
$\mathrm{IC}_2=(-0.6,0.2,-0.2,-0.1,0.1,-0.2)$ (bottom). Solid lines: closed-loop
feedback from the regressed value function; dash-dot lines: open-loop optimum for horizon T
(Pontryagin BVP).}
\label{fig:sat6d_panel}
\end{figure}

At $N=70$ with basis $\mathfrak{I}_{\mathrm{HC}}(4)$ ($q=537$), the zeroth- 
and first-order systems are underdetermined ($M/q\approx 0.13$ and $0.91$ 
respectively), while partial second-order regression with $\rho=0.5$ yields 
$M/q\approx 2.3$, sufficient to recover a stabilising feedback. 
Figure~\ref{fig:sat6d_panel} shows closed-loop trajectories from two initial 
conditions inside $\Omega$; both reach the origin well within the horizon $T=20$.

\subsection{19-dimensional controlled Allen--Cahn equation}
\label{sec:allen_cahn_19d}

We consider a PDE-constrained optimal control problem for the Allen--Cahn 
equation on $\Gamma=(-1,1)$:
\begin{equation}\label{eq:ac_cost}
\min_{\mathbf{u}\in L^2([0,T];\mathbb{R}^m)}
\int_0^T \|y(\cdot,t)\|_{L^2(\Gamma)}^2 
+ \beta\|\mathbf{u}(t)\|_2^2\,dt,
\end{equation}
subject to
\begin{equation}\label{eq:ac_pde_xt}
\begin{cases}
\partial_t y - \nu\partial_x^2 y - y(1-y^2) 
= \displaystyle\sum_{i=1}^m u_i(t)\mathbf{1}_{\omega_i}(x)
& \text{in }\Gamma\times(0,T),\\
\partial_x y(\pm1,t)=0 & \text{in }(0,T),\\
y(x,0)=\tilde{y}_0(x) & \text{in }\Gamma.
\end{cases}
\end{equation}
We discretise in space on a uniform grid of $M-1$ interior nodes with spacing 
$\Delta x=2/M$, denoting by $\bar{y}(t)\in\mathbb{R}^{M-1}$ the nodal values. 
The discretised problem is
\begin{equation}\label{eq:ac_cost_disc}
\min_{\mathbf{u}\in L^2([0,T];\mathbb{R}^m)}
\int_0^T\bar{y}(t)^\top Q_M\bar{y}(t)+\beta\|\mathbf{u}(t)\|_2^2\,dt,
\end{equation}
subject to
\begin{equation}\label{eq:ac_dyn_disc}
\begin{cases}
\partial_t\bar{y} = A_M\bar{y}-\bar{y}\odot\bar{y}\odot\bar{y}+B_M\mathbf{u},\\
\bar{y}_i(0) = \tilde{y}_0(x_i), \quad i=1,\dots,M-1,
\end{cases}
\end{equation}
where $\mathbf{x}:=(\tilde{y}_0(x_1),\dots,\tilde{y}_0(x_{M-1}))\in\Omega
\subset\mathbb{R}^{M-1}$ is the initial condition and
\begin{align*}
A_M &:= \frac{\nu}{\Delta x^2}\bigl(\mathrm{tridiag}(1,-2,1)
+e_1e_1^\top+e_{M-1}e_{M-1}^\top\bigr)+I_{M-1},\\
B_M &:= \bigl(\mathbf{1}_{\omega_1}(\bar{x})\;\cdots\;
\mathbf{1}_{\omega_m}(\bar{x})\bigr),\qquad
Q_M := \Delta x\,I_{M-1}.
\end{align*}
The rank-one corrections in $A_M$ encode the Neumann boundary conditions. 
The PMP adjoint system is
\begin{equation}\label{eq:ac_adjoint}
\begin{cases}
\partial_t p = -\bigl(A_M-3\operatorname{diag}(\bar{y}\odot\bar{y})\bigr)^\top p
-2Q_M\bar{y}, \\
p(T)=0,
\end{cases}
\end{equation}
with optimal control $\mathbf{u}^*=-(2\beta)^{-1}B_M^\top p$. The Riccati 
coefficient matrices \eqref{eq:riccati} are
\begin{equation}\label{eq:ac_hamiltonian_blocks}
\begin{array}{ll}
\mathcal{H}_{py} = A_M-3\operatorname{diag}(\bar{y}\odot\bar{y}), &
\mathcal{H}_{yp} = A_M-3\operatorname{diag}(\bar{y}\odot\bar{y}),\\[4pt]
\mathcal{H}_{pp} = -\dfrac{1}{2\beta}B_MB_M^\top, &
\mathcal{H}_{yy} = -6\operatorname{diag}(p\odot\bar{y})+2\Delta x\,I_{M-1}.
\end{array}
\end{equation}

We set $M=20$, giving $n=19$, and fix $m=3$, $T=4$, $\nu=0.1$, $\beta=0.01$, 
with actuator supports $\omega_1=(-0.7,-0.4)$, $\omega_2=(-0.2,0.2)$, 
$\omega_3=(0.4,0.7)$. The value function is approximated on $\Omega=[-1,1]^{19}$ 
using $\mathfrak{I}_{\mathrm{HC}}(4)$ with $q=14{,}213$ Legendre basis functions, 
ridge regularisation $\lambda=10^{-10}$, and weights $\gamma_1=1.0$, 
$\gamma_2=0.3$. Each sample now contributes $210$ equations to the second-order 
system against $20$ to the first-order one; the zeroth-order system, at one 
equation per sample, is severely underdetermined throughout and is omitted. 
Errors are averaged over $5$ Monte Carlo realisations.

\begin{figure}[h!]
\centering
\includegraphics[width=\textwidth]{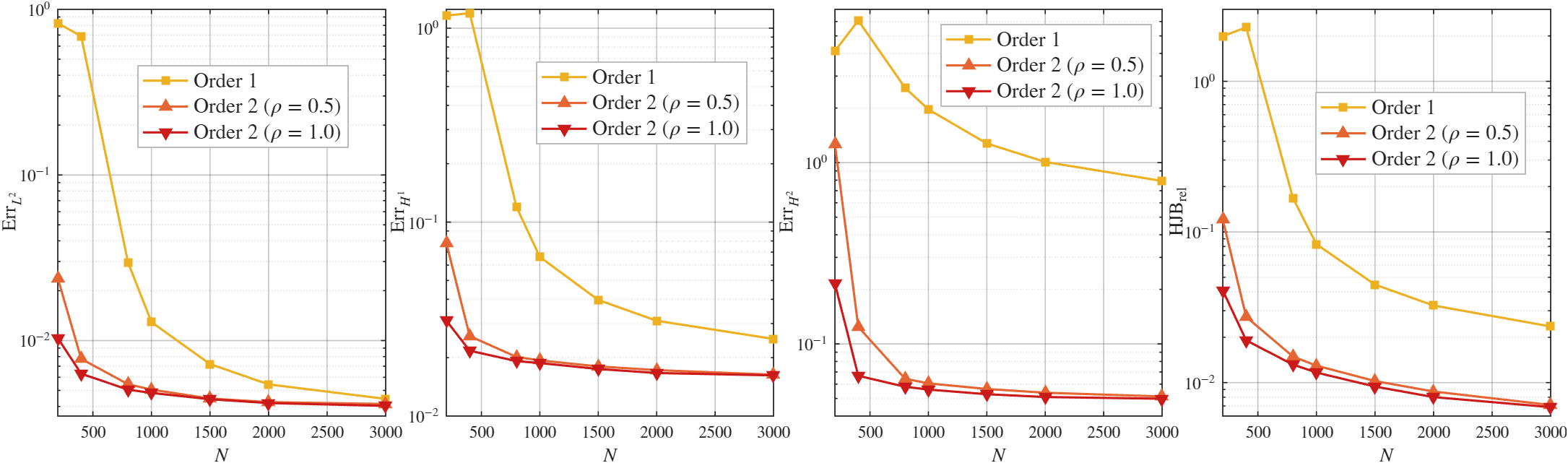}
\caption{Relative validation errors for the $19$-dimensional Allen--Cahn 
system on $\Omega=[-1,1]^{19}$ as a function of $N$. First-order regression 
and second-order regression with $\rho\in\{0.5,1.0\}$ are compared, averaged 
over $5$ Monte Carlo realisations. Basis $\mathfrak{I}_{\mathrm{HC}}(4)$, 
$q=14{,}213$, $\lambda=10^{-10}$. Note that each panel reports a different 
error metric, identified by the subscript of $\mathrm{Err}$, and that the 
ordinate scales differ accordingly across panels.}
\label{fig:19d_final}
\end{figure}

Figure~\ref{fig:19d_final} shows that the $n(n+1)/2=190$ additional equations 
per sample provided by the Hessian block are decisive in this regime. At 
$N=200$, first-order regression yields $\mathrm{Err}_{L^2}=8.3\times10^{-1}$, 
nearly two orders of magnitude above the second-order result of 
$1.0\times10^{-2}$; the second-order system is already well-determined 
($M/q\approx3$) while the first-order system remains underdetermined. As $N$ 
grows, the first-order curves approach the second-order ones in 
$\mathrm{Err}_{L^2}$, $\mathrm{Err}_{H^1}$, and $\mathrm{HJB}_{\mathrm{rel}}$, 
but require substantially larger budgets to do so. The partial Hessian 
configuration $\rho=0.5$ matches the full second-order formulation to within 
$3\%$ at $N=3000$, halving the size of the second-order block at negligible 
accuracy cost.

What distinguishes the $n=19$ regime from lower dimensions is that 
second-order augmentation is no longer merely the most accurate option: at 
any computationally accessible training budget, it is the only formulation 
that brings the regression system into an overdetermined regime on a basis 
expressive enough to approximate $V$. The $n(n+1)/2$ amplification of 
equations per sample thus shifts from being an accuracy advantage to being 
a structural necessity against the curse of dimensionality.
\noindent\subsection*{Concluding remarks}\label{sec:conclusion}

We have proposed a data-driven methodology for value-function and feedback-law 
approximation in deterministic optimal control problems governed by 
Hamilton--Jacobi--Bellman equations. The approach exploits the link between 
the Pontryagin Maximum Principle and dynamic programming to generate augmented 
datasets containing values, gradients, and Hessians of the value function along 
optimal trajectories, processed via sparse polynomial regression with hyperbolic 
cross bases and ridge regularisation.

The numerical experiments highlight the distinctive role of second-order 
information. Beyond refining lower-order regression, Hessian augmentation 
populates approximation spaces unreachable to lower-order formulations at a 
given training budget, with equation-count gains scaling as $\sim n/2$ over 
first-order and $\sim n^2/2$ over zeroth-order data. This advantage sharpens 
with dimension: in the $19$D Allen--Cahn setting, where lower-order systems 
remain underdetermined at moderate $N$, second-order augmentation was the only 
formulation producing a meaningful approximation. The most accurate results 
arose at moderate overdetermination of the second-order system, a balance 
that becomes increasingly accessible as $n$ grows precisely because of the 
$n(n+1)/2$ amplification factor.

The method is causality-free and exposes every ingredient as a design variable: 
sampling, basis, weights, and regularisation. The dominant offline cost is 
dataset generation; online feedback evaluation is linear in $q$ with no 
automatic differentiation required.

Three natural extensions suggest themselves: a theoretical analysis of 
convergence rates and sample complexity for second-order augmented regression; 
adaptive strategies for the partial Hessian formulation by selecting samples by local geometry of 
the value function rather than uniformly; and extensions to stochastic and 
mean-field control problems, where second-order structure is intrinsic to 
the associated HJB equations.

\noindent\subsection*{Acknowledgments}This material is based upon work supported by the Air Force Office of Scientific Research under award number FA8655-26-1-B011.

\renewcommand{\refname}{References}

\end{document}